\documentclass[12pt,oneside]{amsart}
\newtheorem{thm}{Theorem}[section] %uncomment to associate Theorem number to section, e.g. Theorem 1.2
\newtheorem{cor}[thm]{Corollary}
\newtheorem{lem}[thm]{Lemma} % {Lemma}[thm] for numbering within each subsection!!!! i.e. 2.2.1
\newtheorem{prop}[thm]{Proposition}
\theoremstyle{definition}
\theoremstyle{remark}
\newtheorem{rem}[thm]{Remark} %[section]
\theoremstyle{proof}

\numberwithin{equation}{section}

\newcommand{\norm}[1]{\left\Vert#1\right\Vert}
\newcommand{\abs}[1]{\left\vert#1\right\vert}
\newcommand{\set}[1]{\left\{#1\right\}}
\newcommand{\brac}[1]{\left(#1\right)}
\newcommand{\scalar}[1]{\left \langle #1 \right \rangle}
\newcommand{\sscalar}[1]{\langle #1 \rangle}
\newcommand{\Real}{\mathbb{R}}

\newcommand{\To}{\longrightarrow}

\newcommand{\BP}{\mathcal{BP}}
\newcommand{\I}{\mathcal{I}}
\newcommand{\E}{\mathcal{E}}
\newcommand{\M}{\mathcal{M}}
\newcommand{\C}{\mathcal{C}}
\newcommand{\Vol}[1]{\textrm{Vol} (#1)}

\renewcommand{\S}{\mathcal{S}}

\begin{document}

\title{Generalized Intersection Bodies}
\author{Emanuel Milman\\ \bigskip \fontfamily{cmr} \fontseries{m} \fontshape{sc} \fontsize{10}{0} \selectfont
Department of Mathematics, \\
The Weizmann Institute of Science, \\Rehovot 76100, Israel. \\
\medskip E-Mail: emanuel\_milman@hotmail.com.}
\thanks{Supported in part by BSF and ISF}
\begin{abstract}
We study the structures of two types of generalizations of
intersection-bodies and the problem of whether they are in fact
equivalent.
%in the hope of gaining motivation for their equivalence.
Intersection-bodies were introduced by Lutwak and played a key
role in the solution of the Busemann-Petty problem. A natural
geometric generalization of this problem considered by Zhang, led
him to introduce one type of generalized intersection-bodies. A
second type was introduced by Koldobsky, who studied a different
analytic generalization of this problem. Koldobsky also studied
the connection between these two types of bodies, and noted that
an equivalence between these two notions would completely settle
the unresolved cases in the generalized Busemann-Petty problem.
We show that these classes share many identical structural
properties, proving the same results using Integral Geometry
techniques for Zhang's class and Fourier transform techniques for
Koldobsky's class. Using a Functional Analytic approach, we give
several surprising equivalent formulations for the equivalence
problem, which reveal a deep connection to several fundamental
problems in the Integral Geometry of the Grassmann Manifold.
\end{abstract}

% ----------------------------------------------------------------

\maketitle

%%%%%%%%%%%%%%%%%%%%%%%%%%%%%%%%%%%%%%%%%%%%%%%%%%%%%%%%%%%%%%%%%%%%%%%%%
%%%%%%%%%%%%%%%%%%%%%%%% SECTION: INTRODUCTION %%%%%%%%%%%%%%%%%%%%%%%%%%
%%%%%%%%%%%%%%%%%%%%%%%%%%%%%%%%%%%%%%%%%%%%%%%%%%%%%%%%%%%%%%%%%%%%%%%%%

\section{Introduction}

Let $\Vol{L}$ denote the Lebesgue measure of a set $L \subset
\Real^n$ in its affine hull, and let $G(n,k)$ denote the Grassmann
manifold of $k$ dimensional subspaces of $\Real^n$. Let $D_n$
denote the Euclidean unit ball, and $S^{n-1}$ the Euclidean
sphere. All of the bodies considered in this note will be assumed
to be centrally-symmetric star-bodies, defined by a continuous
radial function $\rho_K(\theta) = \max \{r > 0 \; | \; r\theta \in
K \}$ for $\theta \in S^{n-1}$ and a star-body $K$. We shall deal
with two generalizations of the notion of an \emph{intersection
body}, first introduced by Lutwak in
\cite{Lutwak-dual-mixed-volumes} (see also
\cite{Lutwak-intersection-bodies}). A star-body $K$ is said to be
an intersection body of a star-body $L$, if $\rho_K(\theta) =
\Vol{L \cap \theta^\perp}$ for every $\theta \in S^{n-1}$, where
$\theta^\perp$ is the hyperplane perpendicular to $\theta$.
%$\Vol{K \cap H^\perp} = \Vol{L \cap H}$ for every hyperplane $H \in G(n,n-1)$.
$K$ is said to be an intersection body, if it is the limit in the
radial metric $d_r$ of intersection bodies $\{K_i\}$ of
star-bodies $\{L_i\}$, where $d_r(K_1,K_2) = \sup_{\theta \in
S^{n-1}} \abs{\rho_{K_1}(\theta)- \rho_{K_2}(\theta)}$. This is
equivalent (e.g. \cite{Lutwak-intersection-bodies},
\cite{Gardner-BP-5dim}) to $\rho_K = R^*(d\mu)$, where $\mu$ is a
non-negative Borel measure on $S^{n-1}$, $R^*$ is the dual
transform (as in (\ref{eq:duality111})) to the Spherical Radon
Transform $R:C(S^{n-1}) \rightarrow C(S^{n-1})$, which is defined
for $f\in C(S^{n-1})$ as:
\begin{equation} \label{eq:Radon111}
R(f)(\theta) = \int_{S^{n-1} \cap \theta^\perp} f(\xi)
d\sigma_{\theta}(\xi),
\end{equation}
where $\sigma_{\theta}$ the Haar probability measure on $S^{n-1}
\cap \theta^\perp$.

%This is equivalent (e.g. []) to $\rho_K = R(d\mu)$, where $\mu$
%is a non-negative Borel measure on $S^{n-1}$ and $R$ is the
%Spherical Radon Transform: $R(d\mu)(\theta) = \int_{S^{n-1} \cap
%\theta^\perp} d\mu(\xi)$.

The notion of an intersection body has been shown to be
fundamentally connected to the Busemann-Petty Problem (first posed
in \cite{Busemann-Petty}), which asks whether two
centrally-symmetric convex bodies $K$ and $L$ in $\Real^n$
satisfying:
\begin{equation}\label{eq:Busemann-Petty}
\Vol{K\cap H} \leq \Vol{L \cap H} \; \; \forall H \in G(n,n-1)
\end{equation}
necessarily satisfy $\Vol{K} \leq \Vol{L}$. It was shown in
\cite{Lutwak-intersection-bodies}, \cite{Gardner-BP-5dim} that
the answer is equivalent to whether all convex bodies in $\Real^n$
are intersection bodies, and in a series of results
(\cite{Larman-Rogers}, \cite{Ball-BP}, \cite{Bourgain-BP},
\cite{Giannopoulos-BP}, \cite{Papadimitrakis-BP},
\cite{Gardner-BP-5dim}, \cite{Gardner-BP-3dim},
\cite{Koldobsky-lp-intersection-bodies},
\cite{Zhang-Correction-4dim}, \cite{GKS}) that this is true for
$n\leq 4$, but false for $n \geq 5$.

In \cite{Zhang-Gen-BP}, Zhang considered a generalization of the
Busemann-Petty problem, in which $G(n,n-1)$ in
(\ref{eq:Busemann-Petty}) is replaced by $G(n,n-k)$, where $k$ is
some integer between $1$ and $n-1$. Zhang showed that the
\emph{generalized $k$-codimensional Busemann-Petty problem} is
also naturally associated to another class of bodies, which will
be referred to as \emph{$k$-Busemann-Petty bodies} (note that
these bodies are referred to as \emph{$n-k$-intersection bodies}
in \cite{Zhang-Gen-BP} and \emph{generalized $k$-intersection
bodies} in \cite{Koldobsky-I-equal-BP}), and that the generalized
$k$-codimensional problem is equivalent to whether all convex
bodies in $\Real^n$ are $k$-Busemann-Petty bodies. It was shown in
\cite{Bourgain-Zhang} (see also a correction in
\cite{Rubin-Zhang}), and later in \cite{Koldobsky-I-equal-BP},
that the answer is negative for $k < n-3$, but the cases $k=n-3$
and $k=n-2$ still remain open (the case $k=n-1$ is obviously
true). Several partial answers to these cases are known. It was
shown in \cite{Zhang-Gen-BP} (see also \cite{Rubin-Zhang}) that
when $K$ is a centrally-symmetric convex body of revolution then
the answer is positive for the pair $K,L$ with $k=n-2,n-3$ and any
star-body $L$. When $k=n-2$, it was shown in \cite{Bourgain-Zhang}
that the answer is positive if $L$ is a Euclidean ball and $K$ is
convex and sufficiently close to $L$. Several other
generalizations of the Busemann-Petty problem were treated in
\cite{Rubin-Zhang}, \cite{Zvavitch-BP-Arbitrary-Measures},
\cite{BP-Hyperbolic}, \cite{BP-Hyperbolic-low-dim}.

Before defining the class of $k$-Busemann-Petty bodies we shall
need to introduce the $m$-dimensional Spherical Radon Transform,
acting on spaces of continuous functions as follows:
\begin{eqnarray}
\nonumber R_m: C(S^{n-1}) \To C(G(n,m)) \\
\nonumber R_m(f) (E) = \int_{S^{n-1}\cap E} f(\theta)
d\sigma_E(\theta) ,
\end{eqnarray}
where $\sigma_E$ is the Haar probability measure on $S^{n-1} \cap
E$. It is well known (e.g. \cite{Helgason-Book}) that as an
operator on \emph{even} continuous functions, $R_m$ is injective.
The dual transform is defined on spaces of \emph{signed} Borel
measures $\M$ by:
\begin{eqnarray}
\label{eq:duality111} & R_m^*: \M(G(n,m)) \To \M(S^{n-1}) & \\
\nonumber & \int_{S^{n-1}} f R_m^*(d\mu) = \int_{G(n,m)} R_m(f)
d\mu & \forall f \in C(S^{n-1}),
\end{eqnarray}
and for a measure $\mu$ with continuous density $g$, the
transform may be explicitly written in terms of $g$ (see
\cite{Zhang-Gen-BP}):
\begin{eqnarray}
\nonumber R_m^* g (\theta) = \int_{\theta \in E \in G(n,m)} g(E)
d\nu_{m,\theta}(E) ,
\end{eqnarray}
where $\nu_{m,\theta}$ is the Haar probability measure on the
homogeneous space $\set{E \in G(n,m) \;|\; \theta \in E}$.

%The dual transform is given by (see []):
%\begin{eqnarray}
%\nonumber R_m^*: C(G(n,m)) \To C(S^{n-1}) \;\;\;\; R_m^* g
%(\theta) = \int_{\theta \in H \in G(n,m)} g(H) dH ,
%\end{eqnarray}
%where $dH$ is the Haar probability measure on $G(n-1,m-1)$. Both
%of these transforms may be naturally extended to act on the
%measure spaces which are the duals of $C(S^{n-1})$ and
%$C(G(n,m))$.

We shall say that a body $K$ is a $k$-Busemann-Petty body if
$\rho_K^k = R_{n-k}^*(d\mu)$ as measures in $\M(S^{n-1})$, where
$\mu$ is a non-negative Borel measure on $G(n,n-k)$. We shall
denote the class of such bodies by $\BP_k^n$. Choosing $k=1$, for
which $G(n,n-1)$ is isometric to $S^{n-1} / Z_2$ by mapping $H$
to $S^{n-1}\cap H^\perp$, and noticing that $R$ is equivalent to
$R_{n-1}$ under this map, we see that $\BP_1^n$ is exactly the
class of intersection bodies.

Another generalization of the notion of an intersection body,
which was considered by Koldobsky in \cite{Koldobsky-I-equal-BP},
is that of a \emph{$k$-intersection body}. A star-body $K$ is said
to be a $k$-intersection body of a star-body $L$, if $\Vol{K \cap
H^\perp} = \Vol{L \cap H}$ for every $H \in G(n,n-k)$. $K$ is
said to be a $k$-intersection body, if it is the limit in the
radial metric of $k$-intersection bodies $\{K_i\}$ of star-bodies
$\{L_i\}$. We shall denote the class of such bodies by $\I_k^n$.
Again, choosing $k=1$, we see that $\I_1^n$ is exactly the class
of intersection bodies.

In \cite{Koldobsky-I-equal-BP}, Koldobsky considered the
relationship between these two types of generalizations, $\BP_k^n$
and $\I_k^n$, and proved that $\BP_k^n \subset \I_k^n$ (hence our
reluctance to use the term "generalized $n-k$-intersection
bodies" for $\BP_k^n$). Koldobsky also asked whether the opposite
inclusion is equally true for all $k$ between 2 and $n-2$ (for 1
and $n-1$ this is true). If this were true, as remarked by
Koldobsky, a positive answer to the generalized $k$-codimensional
Busemann-Petty problem for $k \geq n-3$ would follow, since for
those values of $k$ any centrally-symmetric convex body in
$\Real^n$ is known to be a $k$-intersection body
(\cite{Koldobsky-correlation-inequality},\cite{Koldobsky-convex-is-n-3-intersection},
\cite{Koldobsky-I-equal-BP}).

Our first remark in this note is that
%In this note we that
the two classes $\BP_k^n$ and $\I_k^n$ share
many identical structural properties, suggesting that it is
indeed reasonable to believe that $\BP_k^n = \I_k^n$.
%The main purpose of this note is to provide new motivation for why it is
%reasonable to believe that $\BP_k^n = \I_k^n$, by giving several
%equivalent formulations for this question and describing a common
%structure which is shared by these classes.
Some previously known characterizations of these classes and
associated tools are outlined in Section \ref{sec:1}, providing
some intuitive motivation and common ground to start from. Some
of these previously known results are also given simplified
proofs in this section. It turns out that the natural language
for handling the class $\I_k^n$ is the language of Fourier
Transforms of homogeneous distributions, developed extensively by
Koldobsky, while the natural language for the class $\BP_k^n$ is
the language of Integral Geometry and Radon Transforms. In Section
\ref{sec:2} we show that both classes share
% Version 2
%many identical structure properties,
a common structure,
% Version 2
by proving the same results for $\BP_k^n$ (using Grassmann
Geometry techniques) and for $\I_k^n$ (using Fourier Transform
techniques). We define the \emph{$k$-radial sum} of two
star-bodies $L_1,L_2$ as the star-body $L$ satisfying $\rho_L^{k}
= \rho_{L_1}^{k} + \rho_{L_2}^{k}$. For each of these classes
$\C_k^n$, where $\C = \I$ or $\C = \BP$ and $k,l=1,\ldots,n-1$,
we show the following:
%in the corresponding section the following:

\smallskip
\noindent\textbf{Structure Theorem. }\emph{
\begin{enumerate}
\item
$\C_k^n$ is closed under full-rank linear transformations,
$k$-radial sums and taking limit in the radial metric.
\item
$\C_1^n$ is the class of intersection-bodies in $\Real^n$, and
$\C_{n-1}^n$ is the class of all symmetric star-bodies in
$\Real^n$.
\item
Let $K_1 \in \C_{k_1}^n$, $K_2 \in \C_{k_2}^n$ and $l = k_1 + k_2
\leq n-1$. Then the star-body $L$ defined by $\rho_L^{l} =
\rho_{K_1}^{k_1} \rho_{K_2}^{k_2}$ satisfies $L \in \C_{l}^n$. As
corollaries:
\begin{enumerate}
\item
$\C_{k_1}^n \cap \C_{k_2}^n \subset \C_{k_1+k_2}^n$ if $k_1+k_2
\leq n-1$.
\item
$\C_{k}^n \subset \C_{l}^n$ if $k$ divides $l$.
\item
If $K \in \C_{k}^n$ then the star-body $L$ defined by $\rho_L =
\rho_K^{k/l}$ satisfies $L \in \C_{l}^n$ for $l \geq k$.
\end{enumerate}
\item
If $K \in \C_k^n$ then any $m$-dimensional central section $L$ of
$K$ (for $m>k$) satisfies $L \in \C_k^m$.
\end{enumerate}
}

(1) and (2) above are well known and basically follow from the
definitions (or from the characterizations in Section
\ref{sec:1}), but we mention them here for completeness. It
should also be clear that (3) implies the three corollaries
following it: (3a) by using $K_1 = K_2$, (3b) by successively
applying (3a), and (3c) by using $K_2 = D_n$. (3) for $\I_k^n$ was
also noticed independently by Koldobsky, but never published. For
$\BP_k^n$, (4) and (3b) for $k=1$ were proved by Grinberg and
Zhang in \cite{Grinberg-Zhang}. In the same paper, a very useful
characterization of the class $\BP_k^n$ was given (see Section
\ref{sec:1}). Combining it with (3) and (3c), we get as a
corollary the following non-trivial result, which is of
independent interest:

\smallskip
\noindent\textbf{Ellipsoid Corollary. }\emph{For any $1 \leq k
\leq n-1$ and $k$ ellipsoids $\set{\E_i}_{i=1}^k$ in $\Real^n$,
define the body $L$ by:
\[
\rho_L =  \rho_{\E_1} \cdot \ldots \cdot \rho_{\E_k},
\]
and let $k \leq l\leq n-1$. Then there exists a sequence of
star-bodies $\set{L_i}$ which tends to $L$ in the radial metric
and satisfies:
\[
\rho_{L_i} = \rho_{\E^i_1}^{l} + \ldots + \rho_{\E^i_{m_i}}^{l},
\]
where $\set{\E^i_j}$ are ellipsoids.}

\smallskip
Naturally, the case $\E_1 = \ldots = \E_k$ is of particular
interest. In the same spirit, we give a strengthened version of
Grinberg and Zhang's characterization of $\BP_k^n$ in Section
\ref{sec:2}. We remark that (3) from the Structure Theorem may in
fact be a characterization of the classes $\I_k^n$ or $\BP_k^n$
for $k>1$. In other words, it may be that for $\C = \BP$ or $\C =
\I$, $L \in \C_k^n$ iff there exist $\set{K_i}_{i=1}^k \subset
\C_1^n$, such that $\rho_L^k = \rho_{K_1} \cdot \ldots \cdot
\rho_{K_k}$. Since in either case $\C_1^n$ is the class of
intersection bodies in $\Real^n$, a proof of such a
characterization for $\C = \I$ and a fixed $k$ would imply that
$\BP_k^n = \I_k^n$ for that $k$.

In order to prove (3) for $\C = \BP$, we derive (what seems to
be) a new formula for integration on products of Grassmann
manifolds. The complete formulation and proof are given in the
Appendix. A very similar formulation of the case
$k_1,\ldots,k_r=1$ was given by Blashcke and Petkantschin (see
\cite{Santalo-Book},\cite{Miles} for an easy derivation), and
used by Grinberg and Zhang in \cite{Grinberg-Zhang} to deduce
that $\BP_1^n \subset \BP_l^n$ for all $1\leq l \leq n-1$. For $F
\in G(n,n-l)$ and $1 \leq k<l \leq n-1$, we denote by
$G_F(n,n-k)$ the manifold $\set{E \in G(n,n-k) | F \subset E}$.
The volume of the parallelepiped mentioned in the statement below
is defined in the Appendix. A simplified formulation then reads as
follows:

\smallskip
\noindent\textbf{Integration on products of Grassmann manifolds. }
\emph{Let $n>1$ be fixed. For $i=1,\ldots,r$, let $k_i \geq 1$
denote integers whose sum $l$ satisfies $l \leq n-1$. For
$a=1,\ldots,n$ denote by $G^{a} = G(n,n-a)$, and by $\mu^{a}$ the
Haar probability measure on $G^a$. For $F \in G^l$ and $a =
1,\ldots,l-1$, denote by $\mu^{a}_F$ the Haar probability measure
on $G^{a}_F$. Denote by $\bar{E} = (E_1,\ldots,E_r)$ an ordered
set with $E_i \in G^{k_i}$. Then for any continuous function
$f(\bar{E}) = f(E_1,\ldots,E_r)$ on $G^{k_1} \times \ldots \times
G^{k_r}$:
\begin{eqnarray}
\nonumber & & \!\!\!\!\!\!\!\!\! \int_{E_1 \in G^{k_1}} \cdots
\int_{E_r\in G^{k_r}}
f(\bar{E}) d\mu^{k_1}(E_1)\cdots d\mu^{k_r}(E_r) = \\
\nonumber & & \!\!\!\!\!\!\!\!\! \int_{F \in G^l} \int_{E_1 \in
G^{k_1}_F} \cdots \int_{E_r \in G^{k_r}_F} f(\bar{E})
\Delta(\bar{E}) d\mu^{k_1}_F(E_1) \cdots d\mu^{k_r}_F(E_r)
d\mu^{l}(F),
\end{eqnarray}
where $\Delta(\bar{E}) = C_{n,\set{k_i},l}
\Omega(\bar{E})^{n-l}$, $C_{n,\set{k_i},l}$ is a constant
depending only on $n,\set{k_i},l$, and $\Omega(\bar{E})$ denotes
the $l$-dimensional volume of the parallelepiped spanned by unit
volume elements of $E_1^\perp,\ldots,E_r^\perp$. }

\smallskip

In Section \ref{sec:3} we attempt to bridge the gap between the
the languages of Integral Geometry and Fourier Transforms, by
establishing several new identities. As a by-product, we show,
for instance, that $Ker R_{n-k}^* = Ker (I\circ R_k)^*$, where
$I:C(G(n,k)) \rightarrow C(G(n,n-k))$ denotes the operator
defined as $I(f)(E) = f(E^\perp)$. Essentially using the latter
result, we show the following equivalence:

\smallskip \noindent\textbf{Equivalence between $k$ and $n-k$. }
\[
\BP_k^n = \I_k^n  \; \textrm{  iff  } \; \BP_{n-k}^n = \I_{n-k}^n.
\]

In Section \ref{sec:4} we try to attack the $\BP_k^n = \I_k^n$
question using the results of the previous sections together with
a functional analytic approach. Our results indicate that this
question is deeply connected to several fundamental questions in
Integral Geometry concerning the structure of the Grassmann
manifold.
%on the geometry of the Grassmann manifold.
Let $C_{+}(S^{n-1})$ denote the set of non-negative continuous
functions on the sphere, and let $R_{n-k}(C(S^{n-1}))_{+}$ denote
the set of non-negative functions in the image of $R_{n-k}$. Let
$\overline{A}$ denote the closure of a set $A$ in the
corresponding normed space. If $\mu \in \M(G(n,n-k))$, let
$\mu^\perp \in \M(G(n,k))$ denote the measure defined by
$\mu^\perp(A) = \mu(A^\perp)$ for any Borel set $A \subset
G(n,k)$, where $A^\perp = \set{E^\perp | E \in A}$.

\smallskip

Fixing $n$ and $1\leq k \leq n-1$, the main result of Section
\ref{sec:4} is the following:

\smallskip \noindent\textbf{Equivalence Theorem. } \emph{The
following statements are equivalent: \textbf{
\begin{enumerate}
\item
Equivalence of generalizations of intersection-bodies.
\[\BP_k^n = \I_k^n.\]
\item
Characterization of non-negative range of $R_{n-k}$.
\smallskip
\begin{equation} \label{eq:non-negative-range}
\overline{R_{n-k}(C(S^{n-1}))_{+}} =
\overline{R_{n-k}(C_{+}(S^{n-1})) + I \circ R_k (C_{+}(S^{n-1}))}.
\end{equation}
\item
A Negation Statement. \\
\textmd{There does not exist a non-negative measure $\mu \in
\M(G(n,n-k))$ such that $R_{n-k}^*(d\mu) \geq 1$ and
$R_k^*(d\mu^\perp) \geq 1$ (where ``$\nu \geq 1$" means that
$\nu-1$ is a non-negative measure), and such that:
\[
\inf \set{\scalar{\mu,f} | f \in R_{n-k}(C(S^{n-1}))_{+} \textrm{
and } \scalar{1,f} = 1 } = 0.
\]
}
\end{enumerate}
} }

The approach developed in Section \ref{sec:3} easily shows (once
again) that $\BP_k^n \subset \I_k^n$. Analogously, it will be
evident that the right hand side of (\ref{eq:non-negative-range})
is a subset of the left hand side.

\medskip

We will say that a set $Z \subset G(n,n-k)$ satisfies the
\emph{covering property} if:
%\smallskip
%\noindent\textbf{The Covering Property. }
\begin{equation} \label{eq:condition-on-A}
\bigcup_{E \in Z} E \cap S^{n-1} = S^{n-1} \; \text{ and } \;
\bigcup_{E \in Z} E^\perp \cap S^{n-1} = S^{n-1}.
\end{equation}

The following natural conjecture is given in Section \ref{sec:4}
(see Lemma \ref{lem:covering-necessary} and Remark
\ref{rem:prove-covering-conj}):

\smallskip
\noindent\textbf{Covering Property Conjecture. } \emph{For any
$n>0$, $1\leq k \leq n-1$, if $Z \subset G(n,n-k)$ is a closed set
satisfying $\bigcup_{E \in Z} E \cap S^{n-1} = S^{n-1}$, then
there exists a non-negative measure $\mu \in \M(G(n,n-k))$
supported in $Z$, such that $R_{n-k}^*(d\mu) \geq 1$.}
\smallskip

\noindent Using this conjecture, we extend formulations (1)-(3)
from the Equivalence Theorem in the following:

\medskip \noindent\textbf{Weak Equivalence Theorem. }
\emph{The following statements are equivalent to each
other:\textbf{
\begin{enumerate}
\setcounter{enumi}{3}
\item
``Injectivity" of the Restricted Radon Transform. \\
\textmd{For any $g \in \overline{R_{n-k}(C(S^{n-1}))_{+}}$, if $Z
= g^{-1}(0)$ satisfies the covering property then $g=0$.}\\
\item
Existence of barely balanced measures.\\
\textmd{For any closed $Z \subset G(n,n-k)$ with the covering
property, there exists a measure $\mu \in \M(G(n,n-k))$ such that
$\mu |_{Z^C} \geq 1$ and $R_{n-k}^*(d\mu) = 0$.}
\end{enumerate}
} \noindent Assuming the Covering Property Conjecture,
formulations (1)-(3) imply (4)-(5). }

\medskip

For us, the formulation in (5) seems to have the most potential
for understanding this problem, although we have not been able to
advance in this direction. Without a doubt, (2) is the most
elegant formulation, and perhaps the most natural for Integral
Geometrists.

\medskip

We conclude by proposing another natural problem in Integral
Geometry. Consider the operator $V_k: C(G(n,k)) \rightarrow
C(G(n,k))$ defined as $V_k = I \circ R_{n-k} \circ R_k^*$. It is
easy to see from general principles of Functional Analysis that
$Ker V_k$ is orthogonal to $\overline{Im V_k}$, and therefore as
an operator from $\overline{Im V_k}$ to itself, $V_k$ is
injective and onto a dense set. We show in Section \ref{sec:3}
that in addition, $V_k$ is self-adjoint. In the case $k=1$,
$C(G(n,1))$ may be identified with the class of even continuous
functions on the sphere $C_e(S^{n-1})$, in which case $V_1 :
C_e(S^{n-1}) \rightarrow C_e(S^{n-1})$ becomes the classical
Spherical Radon Transform $R$ given by (\ref{eq:Radon111}).
Elegant inversion formulas for $V_1$ have been developed by many
authors (see \cite{Helgason-Book} and also \cite{Strichartz},
\cite{Grinberg-Radon-Image}, \cite{Grinberg-Rubin},
\cite{Semyanistyi}, \cite{Petrov}). Is it possible to do the same
for the general $V_k$?

\medskip

\noindent \textbf{Acknowledgments.} I would like to deeply thank
my supervisor Prof. Gideon Schechtman for many informative
discussions, carefully reading the manuscript, and especially for
believing in me and allowing me to pursue my interests. I would
also like to thank Prof. Alexander Koldobsky for going over the
manuscript and for his helpful remarks. I also thank Prof. Semyon
Alesker for helpful information and references about Radon
Transforms.

\section{Additional Notations and Previous Results}
\label{sec:1}

In this section we present some previously known results which
will be useful for us later on. For completeness, we try to at
least sketch the proofs of the main results, and on some
occasions, provide alternative proofs. We also add several useful
notations along the way.

\subsection{Additional Notations}

Let $G$ denote any locally compact topological space. The spaces
of continuous and non-negative continuous real-valued functions
on $G$ will be denoted by $C(G)$ and $C_+(G)$, respectively. When
$G$ has a natural involution operator ``$-$", we will denote by
$C_e(G)$ the space of continuous even functions on $G$. Whenever
it makes sense, we will denote by $C^\infty(G)$ the space of
infinitely smooth real-valued functions on $G$, and define
$C_+^\infty(G)$ and $C_{+,e}^\infty(G)$ accordingly. Similarly,
the spaces of signed and non-negative finite Borel measures on
$G$ will be denoted $\M(G)$ and $\M_+(G)$, respectively. When a
natural involution operator ``$-$" exists, the spaces $\M_e(G)$
and $\M_{+,e}(G)$ will denote the corresponding spaces of even
measures. A measure $\mu$ is called even if $\mu(A) = \mu(-A)$
for every Borel set $A \subset G$. For $\mu \in \M(G)$ and $f \in
C(G)$, we denote by $\scalar{\mu,f}_G$ the action of the measure
$\mu$ on $f$ as a linear functional. Whenever it is clear from
the context what the underlying space $G$ is, we will write
$\scalar{\mu,f}$ instead of $\scalar{\mu,f}_{G}$.

We will always assume that a fixed Euclidean structure is given
on $\Real^n$, and denote by $\abs{x}$ the Euclidean norm of $x
\in \Real^n$. We will denote by $O(n)$ the group of orthogonal
rotations in $\Real^n$. The group of volume-preserving linear
transformations in $\Real^n$ will denoted by $SL(n)$. For $T\in
SL(n)$, we denote $T^{-*} = (T^{-1})^*$.

We will always use $\sigma$ to denote the Haar probability
measure on $S^{n-1}$. $G(n,0)$ and $G(n,n)$ will denote the
trivial atomic manifolds, and these are equipped of course with
the trivial Haar probability measure.

For a star-body $K$ (not necessarily convex), we define its
Minkowski functional as $\norm{x}_K = \min \set{ t \geq 0 \; | \;
x \in t K}$. When $K$ is a centrally-symmetric convex body, this
of course coincides with the natural norm associated with it.
Obviously $\rho_K(\theta) = \norm{\theta}^{-1}_K$ for $\theta \in
S^{n-1}$.

\subsection{Closure under basic operations}

It is not hard to check \emph{from the definitions} that the
classes $\BP_k^n$ and $\I_k^n$ are closed under $k$-radial sums,
full-rank linear transformations and limit in the radial metric.
Indeed, the closure under limit in the radial metric follows from
the definition of $\I_k^n$ and from the $w^*$-compactness of the
unit ball of $\M(G(n,n-k))$ for $\BP_k^n$. The closure under
$k$-radial sums is also immediate for $\BP_k^n$, but for $\I_k^n$
this requires a little more thought. Indeed, by polar
integration, if $K_i$ is a $k$-intersection body of a star-body
$L_i$, for $i=1,2$, then the body $K$ which is the $k$-radial sum
of $K_1$ and $K_2$ is a $k$-intersection body of the $n-k$-radial
sum of $L_1$ and $L_2$, and the general case follows by passing
to a limit. The closure under full-rank linear-transformations
requires a little more ingenuity. It is not so hard to check that
if $K$ is a $k$-intersection body of a star-body $L$ then $T(K)$
is a $k$-intersection body of $T^{-*}(L)$ for $T\in SL(n)$, which
settles the case of $\I_k^n$. For $\BP_k^n$, this requires
additional work, and is actually a good exercise to show
directly. Instead, we prefer to trivially deduce this from
Theorem \ref{thm:G&Z} below.

\subsection{The class $\BP_k^n$} \label{sec:previous-BP}

%but this is an immediate consequence of the
The following characterization of $\BP_k^n$, first proved by
Goodey and Weil in \cite{Goodey-Weil} for intersection-bodies
(the case $k=1$), and extended to general $k$ by Grinberg and
Zhang in \cite{Grinberg-Zhang}, is extremely useful:

\begin{thm}[Grinberg and Zhang] \label{thm:G&Z}
A star-body $K$ is a $k$-Busemann-Petty body iff it is the limit
of $\set{K_i}$ in the radial metric, where each $K_i$ is a finite
$k$-radial sums of ellipsoids $\set{\E^i_j}$:
\[
\rho^k_{K_i} = \rho^k_{\E^i_1} + \ldots + \rho^k_{\E^i_{m_i}}.
\]
\end{thm}

Before commenting on the proof of this theorem, we introduce the
following useful notion used by Grinberg and Zhang. For any $G$, a
homogeneous space of $O(n)$, and measures $\mu \in \M(G)$ and
$\eta \in \M(O(n))$, we define their convolution $\eta \ast \mu
\in \M(G)$ as the measure satisfying $\eta \ast \mu (A) =
\int_{O(n)} \mu(u^{-1}(A)) d\eta(u)$ for every Borel subset $A
\subset G$. The definition is essentially the same when $\eta \in
\M(H)$, where $H$ is another homogeneous space of $O(n)$, by
identifying between $\eta$ and its lifting $\tilde{\eta} \in
\M(O(n))$ defined as $\tilde{\eta}(A) = \eta(\pi(A))$ for any
Borel subset $A \subset O(n)$, where $\pi: O(n) \rightarrow H$ is
the canonical projection.

Let $\sigma_F$ denote the Haar probability measure on
$S^{n-1}\cap F$, so that as a linear functional, for any $f \in
C(S^{n-1})$, $\sigma_F(f) = R_{n-k}(f)(F)$. The key idea
underlying Theorem \ref{thm:G&Z} is an important observation: for
any $F \in G(n,n-k)$, one may explicitly construct a family of
ellipsoids $\set{\E_i(F,\epsilon)}$, such that
$\rho^k_{\E_i(F,\epsilon)}$ tends to $\sigma_F$ in the
$w^*$-topology (as $\epsilon \rightarrow 0$). The ellipsoid
$\E_i(F,\epsilon)$ is defined by:
\[
\norm{x}_{\E_i(F,\epsilon)}^2 =
\frac{\abs{Proj_F(x)}^2}{a(\epsilon)^2} +
\frac{\abs{Proj_{F^\perp}(x)}^2}{b(\epsilon)^2},
\]
where $Proj_E$ denotes the orthogonal projection onto $E$, and
$a(\epsilon),b(\epsilon)$ are chosen appropriately. As observed
by Grinberg and Zhang, one may write $R_{n-k}^* (d\mu)$ as $\mu
\ast \sigma_{F_0}$, where $F_0 = \pi(e)$, $e$ is the identity
element in $O(n)$ and $\pi$ is the canonical projection as above.
Since in the $w^*$-topology, $\sigma_{F_0}$ may be approximated by
$\rho^k_{\E_i(F_0,\epsilon)}$, and $\mu$ by a discrete measure,
the Theorem follows after several technicalities are treated.

We mention a different way to conclude the theorem. It is easy to
verify that:
\[
R_{n-k}(\rho^k_{\E(F,\epsilon)})(E) =
R_{n-k}(\rho^k_{\E(E,\epsilon)})(F) \;\; \forall \; E,F \in
G(n,n-k).
\]
Denoting $G = G(n,n-k)$ for short, if $\rho_K^k =
R_{n-k}^*(d\mu)$ then:
\begin{eqnarray}
\nonumber & & \!\!\!\!\!\!\! R_{n-k}(\rho^k_K)(F) =
\int_{S^{n-1}} \rho^k_K(\theta) d\sigma_F(\theta) =
\lim_{\epsilon \rightarrow 0} \int_{S^{n-1}}
\rho^k_{\E(F,\epsilon)} (\theta) \rho^k_K
(\theta) d\sigma(\theta) =  \\
\nonumber &  & \!\!\!\!\!\!\!\!\! \lim_{\epsilon \rightarrow 0}
\int_{S^{n-1}} \rho^k_{\E(F,\epsilon)} (\theta)
R_{n-k}^*(d\mu)(\theta) d\sigma(\theta) = \lim_{\epsilon
\rightarrow 0} \int_{G}
R_{n-k}(\rho^k_{\E(F,\epsilon)})(E) d\mu(E) =  \\
\nonumber &  &  \lim_{\epsilon \rightarrow 0} \int_{G}
R_{n-k}(\rho^k_{\E(E,\epsilon)})(F) d\mu(E) = R_{n-k}\brac{
\lim_{\epsilon \rightarrow 0} \int_{G} \rho^k_{\E(E,\epsilon)}
d\mu(E) }(F),
\end{eqnarray}
where we have used the uniform convergence of all the limits
involved and that $R_{n-k}$ is a continuous operator w.r.t the
maximum-norm. The result then follows from the injectivity of
$R_{n-k}$ on $C_e(S^{n-1})$.

\medskip Grinberg and Zhang's characterization of the class
$\BP_k^n$ implies that it is actually generated from $D_n$, the
Euclidean unit Ball, by taking full-rank linear transformations,
$k$-radial sums, and limit in the radial metric. By starting from
any other star-body $L$ and performing these operations, it is
obvious that $D_n$ may be constructed, and therefore we see that
$\BP_k^n$ is the minimal non-empty class which is closed under
these three operations. Since $\I_k^n$ trivially contains $D_n$
and is also closed under these operations, it immediately follows
that:
\begin{cor} \label{cor:BP-subset-I}
$\BP_k^n \subset \I_k^n$.
\end{cor}
\noindent This was first observed by Koldobsky in
\cite{Koldobsky-I-equal-BP} using a different approach. We will
give another proof of this in Corollary
\ref{cor:BP-subset-I-again}, which is in a sense more concrete.

\medskip

We conclude this preliminary discussion of the class $\BP_k^n$ by
elaborating a little more on the operation of convolution between
measures on homogeneous spaces of $O(n)$. Let $G,H$ denote
homogeneous spaces of $O(n)$. We identify between a function $f
\in C(G)$ and the measure on $C(G)$ whose density w.r.t. the Haar
probability measure on $G$ is given by $f$, and consider
expressions of the form $f \ast \mu$ and $\mu \ast f$ for $\mu
\in \M(H)$. With the same notations, if $f \in C^\infty(G)$  then
a standard argument shows that $f \ast \mu \in C^\infty(H)$ and
that $\mu \ast f \in C^\infty(G)$. If $\eta \in \M(O(n))$, it is
immediate to check that $\scalar{\mu, \eta \ast f}_G =
\scalar{\eta^{-1} \ast \mu, f}_G$, where $\eta^{-1}\in \M(O(n))$
is the measure defined by $\eta^{-1}(A) = \eta(A^{-1})$ and
$A^{-1} = \set{u^{-1} | u \in A}$ for a Borel set $A \subset
O(n)$. If $\mu_i \in \M(G_i)$ for $i=1,2,3$, one may verify that
this operation is associative: $(\mu_1 \ast \mu_2) \ast \mu_3 =
\mu_1 \ast (\mu_2 \ast \mu_3)$. We conclude with the following
lemma from \cite{Grinberg-Zhang} which will be useful later on.

\begin{lem} \label{lem:approximate-identity}
There exists a sequence of functions $\set{u_i} \subset
C_+^\infty(O(n))$ called an \emph{approximate identity}, such
that for any homogeneous space $G$ of $O(n)$:
\begin{enumerate}
\item
For any $\mu \in \M(G)$, $u_i \ast \mu \in C^\infty(G)$ tends to
$\mu$ in the $w^*$-topology.
\item
For any $g \in C(G)$, $u_i \ast g \in C^\infty(G)$ tends to $g$
uniformly.
\end{enumerate}
\end{lem}

\subsection{The class $\I_k^n$}

In order to handle the class $\I_k^n$, we shall need to adopt a
technique extensively used by Koldobsky: Fourier transforms of
homogeneous distributions. We will only outline the main ideas
here, usually omitting the technical details - we refer the
reader to \cite{Koldobsky-Book} for those. We denote by
$\S(\Real^n)$ the space of rapidly decreasing infinitely
differentiable test functions in $\Real^n$, and by $\S'(\Real^n)$
the space of distributions over $\S(\Real^n)$. The Fourier
Transform $\hat{f}$ of a distribution $f \in \S'(\Real^n)$ is
defined by $\sscalar{\hat{f},\phi} = \sscalar{f,\hat{\phi}}$ for
every test function $\phi$, where $\hat{\phi}(y) = \int \phi(x)
\exp(-i\sscalar{x,y}) dx$. A distribution $f$ is called
homogeneous of degree $p \in \Real$ if $\scalar{f,\phi(\cdot/t)} =
\abs{t}^{n+p} \scalar{f,\phi}$ for every $t>0$, and it is called
even if the same is true for $t=-1$. An even distribution $f$
always satisfies $(\hat{f})^\wedge = (2\pi)^n f$. The Fourier
Transform of an even homogeneous distribution of degree $p$ is an
even homogeneous distribution of degree $-n-p$. A distribution
$f$ is called positive if $\sscalar{f,\phi} \geq 0$ for every
$\phi \geq 0$, implying that $f$ is necessarily a non-negative
Borel measure on $\Real^n$. We use Schwartz's generalization of
Bochner's Theorem (\cite{Gelfand-Shilov}) as a definition, and
call a homogeneous distribution positive-definite if its Fourier
transform is a positive distribution.

Before proceeding, let us give some intuition about how the
Fourier transform of a homogeneous continuous function looks
like. Because of the homogeneity, it is enough to consider a
continuous function on the sphere $f \in C(S^{n-1})$, and take
its homogeneous extension of degree $p \in \Real$, denoted
$E_p(f)$, to the entire $\Real^n$ (formally excluding $\set{0}$ if
$p < 0$). When $p > -n$, the function $E_p(f)$ is locally
integrable, and its action as a distribution on a test function
$\phi$ is simply by integration. Passing to polar coordinates, we
have:
\[
\scalar{E_p(f),\phi} = \int_{S^{n-1}} f(\theta) \int_0^\infty
r^{p+n-1} \phi(r \theta) dr d\theta.
\]
When $p \leq -n$, we can no longer interpret the action of
$E_p(f)$ as an integral. Fortunately, we will mainly be concerned
with Fourier transforms of continuous functions which are
homogeneous of degree $p \in (-n,0)$. This ensures that the
Fourier transform is a homogeneous distribution of degree $-p-n$,
which is in the same range $(-n,0)$. Note that the resulting
distribution need not necessarily be a continuous function on
$\Real^n \setminus \set{0}$, nor even a measure on $\Real^n$
(although this will not occur in our context). We will denote by
$E_p^\wedge(f)$ the Fourier transform of $E_p(f)$. In order to
ensure that $E_p^\wedge(f)$ is a continuous function, we need to
add some smoothness assumptions on $f$ (\cite{Koldobsky-Book}).
We remark that for a continuous function $f \in C(S^{n-1})$,
$E_p^\wedge(f)$ is always continuous for $p \in (-n,n+1]$, and
that for an infinitely smooth $f \in C^\infty(S^{n-1})$,
$E_p^\wedge(f)$ is infinitely smooth for any $p\in(-n,0)$.
Whenever $E_p^\wedge(f)$ is continuous on $\Real^n \setminus
\set{0}$, it is uniquely determined by its value on $S^{n-1}$ (by
homogeneity). In that case, by abuse of notation, we identify
between $E_p^\wedge(f)$ and its restriction to $S^{n-1}$, and in
particular, consider $E_p^\wedge$ as an operator from
$C^\infty(S^{n-1})$ to $C^\infty(S^{n-1})$.

When $f=1$, it is easy to verify that $E_p^\wedge(1)$ is
rotational invariant, so by the homogeneity, it must be a multiple
of $E_{-n-p}(1)$. For a rigorous proof we refer to \cite[p. 192]
{Gelfand-Shilov}, and state this for future reference as:
\begin{lem} \label{lem:E_p(1)}
%\label{eq:E_p(1)}
Fix $n$ and let $p\in (0,n)$. Then:
\[
E_{-p}^\wedge(1) = c(n,p) E_{-n+p}(1) \text{ where  } \; c(n,p) =
\pi^{n/2} 2^{n-p}
\frac{\Gamma\brac{\frac{n-p}{2}}}{\Gamma\brac{\frac{p}{2}}}.
\]
Since $(E_{-p}^\wedge(1))^\wedge = (2\pi)^n E_{-p}(1)$, it is
clear that:\[c(n,p)c(n,n-p) = (2\pi)^n.\]
\end{lem}

\smallskip
The following characterization was given by Koldobsky in
\cite{Koldobsky-I-equal-BP}:

\begin{thm}[Koldobsky] \label{thm:I-char}
%\smallskip
%\noindent\textbf{Theorem (Koldobsky).}
The following are equivalent for a centrally-symmetric star-body
$K$ in $\Real^n$:
\begin{enumerate}
\item $K$ is a $k$-intersection body.
\item $\norm{x}_K^{-k}$ is a positive definite distribution on
$\Real^n$, meaning that its Fourier-transform
$(\norm{\cdot}_K^{-k})^{\wedge}$ is a non-negative Borel measure
on $\Real^n$.
\item
The space $(\Real^n,\norm{\cdot}_K)$ embeds in $L_{-k}$.
%The space $(\Real^n,\norm{\cdot}_K)$ formally ``embeds in $L_{-k}$",
%meaning that it satisfies some characterizing property of
%embedding in $L_p$, continued analytically to the negative value
%$p=-k$, using Fourier-transforms of distributions.
\end{enumerate}
\end{thm}

For completeness, we briefly give the definition of embedding in
$L_{-k}$, although we will not use this later on. Let us denote
the class of centrally-symmetric star bodies $K$ in $\Real^n$ for
which $(\Real^n,\norm{\cdot}_K)$ embeds in $L_p$ by $SL_p^n$. For
$p>0$, it is well known (e.g. \cite{Koldobsky-I-equal-BP}) that $K
\in SL_p^n$ iff:
\begin{equation} \label{eq:SL_p^n}
\norm{x}_K^p = \int_{S^{n-1}} \abs{\scalar{x,\theta}}^p
d\mu_K(\theta),
\end{equation}
for some $\mu_K \in \M_+(S^{n-1})$. Unfortunately, this
characterization breaks down at $p=-1$ since the above integral
no longer converges. However, Koldobsky showed that it is
possible to regularize this integral by using Fourier-transforms
of distributions, and gave the following definition:
$(\Real^n,\norm{\cdot}_K)$ embeds in $L_{-p}$ for $0<p<n$ iff
there exists a measure $\mu_K \in \M_+(S^{n-1})$ such that for any
even test-function $\phi$:
\begin{equation} \label{eq:L_-k-characterization}
\int_{\Real^n} \norm{x}_K^{-p} \phi(x) dx = \int_{S^{n-1}}
\int_0^\infty t^{p-1} \hat{\phi}(t \theta) dt d\mu_K(\theta).
\end{equation}

Let us review the statements of Theorem \ref{thm:I-char}. (2) is
an extremely useful characterization of $k$-intersection bodies,
and immediately implies the closure of $\I_k^n$ under the standard
three operations. Characterization (3) provides additional
motivation for why it is reasonable to believe that $\BP_k^n =
\I_k^n$. For $p\neq 0$, the $p$-norm sum of two bodies $L_1,L_2$
is defined as the body $L$ satisfying $\norm{\cdot}_L^p =
\norm{\cdot}_{L_1}^p + \norm{\cdot}_{L_2}^p$. We will denote by
$D_p^n$, the class of bodies created from $D_n$ by applying
full-rank linear-transformations, $p$-norm sums, and taking the
limit in the radial metric. Using the characterization in
(\ref{eq:SL_p^n}), it is easy to show (e.g. \cite[Theorem 6.13]
{Grinberg-Zhang}) that for $p>0$, the class $SL_p^n$ coincides
with $D_p^n$. Although this characterization breaks down at
$p=-1$, it is still reasonable to expect that the property $SL_p^n
= D_p^n$ should pass over to negative values of $p$ when $SL_p^n$
is (in some sense) extended to this range and becomes $SL_{-k}^n
= \I_k^n$. But by Grinberg and Zhang's characterization (Theorem
\ref{thm:G&Z}), this is exactly satisfied by $k$-Busemann-Petty
bodies: $\BP_k^n = D_{-k}^n$. This suggests that indeed $\BP_k^n
= \I_k^n$.

In addition to the characterization (3) of $\I_k^n$ as the class
of unit-balls of subspaces of \emph{scalar} $L_{-k}$ spaces, a
functional analytic characterization of $\BP_k^n$ as the class of
unit-balls of subspaces of \emph{vector valued} $L_{-k}$ spaces
(in a manner similar to (\ref{eq:L_-k-characterization})), was
given in \cite{Koldobsky-I-equal-BP}. This provides additional
motivation for believing that $\BP_k^n = \I_k^n$, as this would
be an extension to negative values of $p$ of the fact that every
separable vector valued $L_p$ space is isometric to a subspace of
a scalar $L_p$ space and vice-versa.

\medskip

We proceed to explain why (1) and (2) in Theorem \ref{thm:I-char}
are equivalent. To this end, we will need the following Spherical
Parseval identity, due to Koldobsky (\cite{Koldobsky-Book}):

\smallskip
\noindent\textbf{Spherical Parseval (Koldobsky).} \emph{Let $f,g
\in C^\infty_e(S^{n-1})$, and $p\in(0,n)$. Then:
%Extend $f$ and $g$ as
%homogeneous functions on $\Real^n \setminus \set{0}$ of degrees
%$-p$ and $-n+p$, respectively.
\[
\int_{S^{n-1}} E_{-p}^\wedge(f)(\theta) E_{-n+p}^\wedge(g)(\theta)
d\sigma(\theta) = (2 \pi)^n \int_{S^{n-1}} f(\theta) g(\theta)
d\sigma(\theta).
\]
}

\noindent We prefer to present a self-contained proof of this
identity, which seems simpler than the previous approaches in
\cite{Koldobsky-Book}.

\begin{proof}
% Give simple proof of the Spherical Parseval identity.
Let $f = \sum_{k=0}^\infty f_{k}$ and $g = \sum_{k=0}^\infty
g_{k}$ be the canonical decompositions into spherical harmonics,
where $f_k,g_k \in H_k$ and $H_k$ is the space of spherical
harmonics of degree $k$. Since $f$ and $g$ are even, it follows
that $f_{2k+1} = g_{2k+1} = 0$. It is well known
(\cite{Stein-Weiss}) that for $q \in (-n,0)$, the linear operator
$E_q^\wedge : C^\infty(S^{n-1}) \rightarrow C^\infty(S^{n-1})$
decomposes into a direct sum of scalar operators acting on $H_k$.
Indeed, one only needs to check that the $H_k$'s are eigenspaces
of $E_q^\wedge$, and by Schur's Representation Lemma and the fact
that the Fourier transform commutes with the action of the
orthogonal group, it follows that $E_q^\wedge$ must act as a
scalar on these spaces. Denote by $c^{(q)}_k$ the eigenvalue
satisfying $E_q^\wedge(h_k) = c^{(q)}_k h_k$ for any $h_k \in
H_k$. The exact value of $c^{(q)}_k$ is well known (\cite[Theorem
4.1]{Stein-Weiss}), but is irrelevant to our proof. It remains to
notice that since:
\[
E_{-n+p}^\wedge(E_{-p}^\wedge(f))  =
\brac{E_{-p}(f)^\wedge}^\wedge |_{S^{n-1}} = (2 \pi)^n f ,
\]
for any $f \in C^\infty_e(S^{n-1})$, we must have $c^{(-n+p)}_{k}
c^{(-p)}_{k} = (2 \pi)^n$ for all \emph{even} $k$'s. Using the
fact that spherical harmonics of different degrees are orthogonal
to each other in $L_2(S^{n-1})$, and that $f,g,E_{-p}^\wedge(f)$
and $E_{-n+p}^\wedge(g)$ are all in $L_2(S^{n-1})$, we conclude:
\begin{eqnarray}
\nonumber & & \!\!\!\!\!\!\!\!\! \int_{S^{n-1}}
E_{-p}^\wedge(f)(\theta) E_{-n+p}^\wedge(g)(\theta)
d\sigma(\theta) = \int_{S^{n-1}} \sum_{k=0}^\infty c_k^{(-p)}
f_k(\theta)
\sum_{l=0}^\infty c_l^{(-n+p)} g_l(\theta) d\sigma(\theta) \\
\nonumber & = & \int_{S^{n-1}} \sum_{k=0}^\infty c_k^{(-p)}
c_k^{(-n+p)} f_k(\theta) g_k(\theta) d\sigma(\theta) = (2 \pi)^n
\int_{S^{n-1}}
\sum_{k=0}^\infty f_k(\theta) g_k(\theta) d\sigma(\theta) \\
\nonumber & = & (2 \pi)^n \int_{S^{n-1}} \sum_{k=0}^\infty
f_k(\theta) \sum_{l=0}^\infty g_l(\theta) d\sigma(\theta) = (2
\pi)^n \int_{S^{n-1}} f(\theta) g(\theta) d\sigma(\theta).
\end{eqnarray}
\end{proof}

%Since the spherical harmonics form a complete system in
%$L_2(S^{n-1})$,
Note that the above argument actually shows that the Spherical
Parseval identity is also valid when $f,g, E_{-p}^\wedge(f),
E_{-n+p}^\wedge(g) \in L_2(S^{n-1})$.

\begin{rem} \label{rem:spherical-parseval}
Applying the theorem to $g = E_{-p}^\wedge(g')$ for $g'\in
C^\infty_e(S^{n-1})$ and using that $E_{-n+p}^\wedge(g) = (2
\pi)^n g'$, we note that the Spherical Parseval identity has the
following equivalent form, which we will sometimes use:
\[
\int_{S^{n-1}} E_{-p}^\wedge(f)(\theta) g(\theta) d\sigma(\theta)
= \int_{S^{n-1}} f(\theta) E_{-p}^\wedge(g)(\theta)
d\sigma(\theta).
\]
\end{rem}

Another useful result due to Koldobsky, which looks very similar
to the Spherical Parseval identity, is the following:

\begin{thm}[Koldobsky] \label{thm:R_k(wedge)}
%\smallskip
%\noindent\textbf{Theorem (Koldobsky).} \emph{
Let $f \in C^\infty_e(S^{n-1})$, and let \mbox{$k =
1,\ldots,n-1$}. Then for any $H \in G(n,k)$:
\[
\int_{S^{n-1} \cap H^\perp} E_{-k}^\wedge(f)(\theta)
d\sigma_{H^\perp}(\theta) = c(n,k) \int_{S^{n-1} \cap H} f(\theta)
d\sigma_{H}(\theta),
\]
where $c(n,k)$ is the constant from Lemma \ref{lem:E_p(1)}.
%}
\end{thm}

Informally, the latter Theorem may be considered as a special
case of the Spherical Parseval identity, by setting $g =
d\sigma_{H}$ and verifying that in the appropriate sense
$E_{-n+k}^\wedge (d\sigma_{H}) = c(n,k) d\sigma_{H^\perp}$. The
constant in front of the right hand integral is verified by
choosing $f=1$ and using Lemma \ref{lem:E_p(1)}. One way to make
this argument work is to use Grinberg and Zhang's approximation
of $d\sigma_H$ by the functions $\rho^{n-k}_{\E_i}$, which when
written as $\norm{\cdot}^{-n+k}_{\E_i}$ are seen to be already
homogeneous of degree $-n+k$. Computing the Fourier transform is
particularly easy, since $\E_i = T_i(D_n)$, and therefore:
\begin{eqnarray}
\nonumber (\norm{\cdot}^{-n+k}_{T_i(D_n)})^\wedge(x) =
(\norm{T_i^{-1}(\cdot)}^{-n+k}_{D_n})^\wedge(x) = det(T_i)
(\norm{\cdot}^{-n+k}_{D_n})^\wedge(T_i^{*}(x)) = \\
\nonumber det(T_i) d(n,k) \norm{T_i^{*}(x)}^{-k}_{D_n} = det(T_i)
d(n,k) \norm{x}^{-k}_{T_i^{-*}(D_n)}.
\end{eqnarray}
Using Grinberg and Zhang's approximation again, it turns out that
$det(T_i) d(n,k) \rho^{k}_{T_i^{-*}(D_n)}$ tends in the
$w^*$-topology to $c(n,k) d\sigma_{H^\perp}$.

\smallskip

We can now sketch a proof of Koldobsky's Fourier transform
characterization of $k$-intersection bodies. By abuse of
notation, when $(\norm{\cdot}_{K}^{-k})^\wedge$ is continuous, we
will often use $\norm{\cdot}_{K}^{-k}$,
$(\norm{\cdot}_{K}^{-k})^\wedge$ to indicate both locally
integrable functions on $\Real^n$ and continuous functions on
$S^{n-1}$. By definition, an \emph{infinitely smooth} star-body
$K$ which is a $k$-intersection body of a star-body $L$,
satisfies $\Vol{K \cap H^\perp} = \Vol{L \cap H}$ for all $H \in
G(n,n-k)$. Passing to polar coordinates, this is equivalent to:
\[
R_k(\norm{\cdot}_K^{-k})(H^\perp) =
\frac{\Vol{D_{n-k}}}{\Vol{D_k}} R_{n-k}(\norm{\cdot}_L^{-n+k})(H)
\;\; \forall H \in G(n,n-k).
\]
But using Theorem \ref{thm:R_k(wedge)}, we see that:
\[
R_k(\norm{\cdot}_K^{-k})(H^\perp) = c(n,k)^{-1}
R_{n-k}((\norm{\cdot}_K^{-k})^\wedge)(H) \;\; \forall H \in
G(n,n-k).
\]
From the injectivity of $R_{n-k}$ on $C_e(S^{n-1})$, it follows
that:
\[
(\norm{\cdot}_K^{-k})^\wedge = c(n,k)
\frac{\Vol{D_{n-k}}}{\Vol{D_k}} \norm{\cdot}_L^{-n+k}
\]
on $S^{n-1}$, and hence on all $\Real^n$ by homogeneity. We
conclude that $(\norm{\cdot}_K^{-k})^\wedge$ is a non-negative
continuous function on $\Real^n \setminus \set{0}$, and hence
positive as a distribution. For an arbitrary star-body $K$ which
is a $k$-intersection body of a star-body $L$, the same conclusion
holds by approximation ($(\norm{\cdot}_K^{-k})^\wedge$ is still
continuous by the continuity of $\norm{\cdot}_L^{-n+k}$). One may
also invert the argument, proving that for a star-body $K$, if
$(\norm{\cdot}_K^{-k})^\wedge$ is a \emph{continuous} function
which is non-negative, then $K$ is a $k$-intersection body of a
star-body $L$ (defined as above). Taking the limit in the radial
metric, $(\norm{\cdot}_K^{-k})^\wedge$ need not necessarily be a
continuous function for a general $k$-intersection body $K$ which
is the limit of the bodies $\set{K_i}$ (which are
$k$-intersection bodies of star-bodies). Nevertheless, the
non-negative continuous functions
$(\norm{\cdot}_{K_i}^{-k})^\wedge$ must satisfy:
\[
\int_{S^{n-1}} \abs{(\norm{\cdot}_{K_i}^{-k})^\wedge(\theta)}
d\sigma(\theta) = c(n,k) \int_{S^{n-1}} \norm{\theta}_{K_i}^{-k}
d\sigma(\theta),
\]
by the Spherical Parseval identity with $g = 1$ and Lemma
\ref{lem:E_p(1)}, and therefore the integral on the left hand
side is bounded. Using the compactness of the unit-ball of
$\M(S^{n-1})$ in the $w^*$-topology, there must be an accumulation
point of $\{(\norm{\cdot}_{K_i}^{-k})^\wedge\}$, which is a
non-negative Borel measure on $S^{n-1}$. This argument is the
main idea in the proof that for a star-body $K$, $K \in \I_k^n$
iff $(\norm{\cdot}_{K}^{-k})^\wedge$ is a non-negative Borel
measure on $\Real^n$.

\medskip

When $K$ is infinitely smooth, we summarize this in the following
alternative definition for $\I_k^n$, and use it instead of the
original one:

\medskip \noindent\textbf{Alternative Definition of $\I_k^n$.}
\emph{For an infinitely smooth star-body $K$, $K\in \I_k^n$ iff
$(\norm{\cdot}_{K}^{-k})^\wedge \geq 0$ as a $C^\infty$ function
on $S^{n-1}$. }
\smallskip

For a general star-body $K$, we will use Koldobsky's
characterization in the following spherical version, which is an
immediate consequence of the above reasoning (a rigorous proof is
given in \cite[Corollary 3.23]{Koldobsky-Book}):

\begin{prop} \label{prop:I-spherical-char}
For a star-body $K$, $K \in \I_k^n$ iff there exists a
non-negative Borel measure $\mu$ on $S^{n-1}$, such that for any
$f \in C^\infty_e(S^{n-1})$:
\[
\int_{S^{n-1}} f(\theta) \rho_K^k(\theta) d\sigma(\theta) =
\int_{S^{n-1}} E_{-n+k}^\wedge(f)(\theta) d\mu(\theta).
\]
\end{prop}

%%%%%%%%%%%%%%%%%%%%%%%%%%%%%%%%%%%%%%%%%%%%%%%%%%%%%%%%%%%%%%%%%%%%%%%%%%%%%%%%%%%%%%%%%%%%%
%%%%%%%%%%%%%%%%%%%%%%%%%%%%%%%%%%%%%%%%%%%%%%%%%%%%%%%%%%%%%%%%%%%%%%%%%%%%%%%%%%%%%%%%%%%%%
%%%%%%%%%%%%%%%%%%%%%%%%%%%%%%%%%%%%%%%%%%%%%%%%%%%%%%%%%%%%%%%%%%%%%%%%%%%%%%%%%%%%%%%%%%%%%
%%%%%%%%%%%%%%%%%%%%%%%%%%%%%%%%%%%%%%%%%%%%%%%%%%%%%%%%%%%%%%%%%%%%%%%%%%%%%%%%%%%%%%%%%%%%%
%%%%%%%%%%%%%%%%%%%%%%%%%%%%%%%%%%%%%%%%%%%%%%%%%%%%%%%%%%%%%%%%%%%%%%%%%%%%%%%%%%%%%%%%%%%%%
%%%%%%%%%%%%%%%%%%%%%%%%%%%%%%%%%%%%%%%%%%%%%%%%%%%%%%%%%%%%%%%%%%%%%%%%%%%%%%%%%%%%%%%%%%%%%
%%%%%%%%%%%%%%%%%%%%%%%%%%%%%%%%%%%%%%%%%%%%%%%%%%%%%%%%%%%%%%%%%%%%%%%%%%%%%%%%%%%%%%%%%%%%%
%%%%%%%%%%%%%%%%%%%%%%%%%%%%%%%%%%%%%%%%%%%%%%%%%%%%%%%%%%%%%%%%%%%%%%%%%%%%%%%%%%%%%%%%%%%%%

\section{The Identical Structures of $\BP_k^n$ and $\I_k^n$}
\label{sec:2}

In this section we will prove the Structure Theorem, which was
formulated in the Introduction. We will skip over item 1 which
basically follows from the definitions, and was already explained
in detail in Section \ref{sec:1}. Item 2 also follows immediately:
by definition, $\I_1^n = \BP_1^n$ is exactly the class of
intersection bodies in $\Real^n$; any star-body $K$ in $\Real^n$
is an $n-1$-intersection body of a star-body $L$, defined by
$\rho_L(\theta) = 1/2 \Vol{K\cap \theta^\perp}$; and by
definition, $R_1^*$ acts as the identity on $C_e(S^{n-1})$, hence
$\rho_K^{n-1} = R_1^*(\rho_K^{n-1})$ for any star-body $K$,
implying that $K \in \BP_{n-1}^n$. We therefore commence the proof
from item 3. We will prove the Theorem for $\BP_k^n$ and $\I_k^n$
separately, because of the different techniques involved in the
proof.

Before we start, we will need the following useful lemma, which
appears implicitly in \cite{Grinberg-Zhang}. We denote by
$\BP_k^{n,\infty}$ the class of star-bodies $K$ such that
$\rho_K^k = R_{n-k}^*(g)$, where $g \in C_+^\infty(G(n,n-k))$.
Obviously $\BP_k^{n,\infty} \subset \BP_k^n$.

\begin{lem}[\cite{Grinberg-Zhang}] \label{lem:BP-dense}
$\BP_k^{n,\infty}$ is dense in $\BP_k^n$. In particular, the
class of infinitely smooth bodies in $\BP_k^n$ is dense in
$\BP_k^n$.
\end{lem}
\begin{proof}
Let $K \in \BP_k^n$, and assume that $\rho_K^k = R_{n-k}^*(d\mu)$
where $d\mu \in \M_+(G(n,n-k))$. Let $\set{u_i} \subset
C^\infty(O(n))$ be an approximate identity as in Lemma
\ref{lem:approximate-identity}. Let $K_i$ be the star-body for
which $\rho_{K_i}^k = u_i \ast \rho_{K}^k$. Then by Lemma
\ref{lem:approximate-identity}, $\set{K_i}$ is a sequence of
infinitely smooth star-bodies which tend to $K$ in the radial
metric. As in the proof of Theorem \ref{thm:G&Z}, we write
$\rho_K^k = \mu \ast \sigma_{H_0}$, and therefore:
\[
\rho_{K_i}^k = u_i \ast (\mu \ast \sigma_{H_0}) = (u_i \ast \mu)
\ast \sigma_{H_0} = R_{n-k}^*(u_i \ast \mu).
\]
Since $u_i \ast \mu \in C_+^\infty(G(n,n-k)$, this concludes the
proof of the lemma.
\end{proof}

\begin{rem} \label{rem:BP-use-smooth}
By the Lemma and the closure of $\BP_k^n$ (for any $k =
1,\ldots,n-1$) under limit in the radial metric, it is enough to
prove all the remaining items for the classes $\BP_k^{n,\infty}$.
\end{rem}

We will also require the following notations. Given $F \in
G(n,m)$ and $k\geq m$, we denote by $G_F(n,k)$ the manifold
$\set{E \in G(n,k) | F \subset E}$. For $\theta \in S^{n-1}$ we
identify between $\theta$ and the one-dimensional subspace
spanned by it. $G_F(n,k)$ is a homogeneous space of $O(n)$,
therefore there exists a unique Haar probability measure on
$G_F(n,k)$, which is invariant to orthogonal rotations in $O(n)$
which preserve $F$. Thus, if we denote by $\nu_\sigma$ the Haar
probability measure on $G_\sigma(n,m)$ for $\sigma \in S^{n-1}$,
then for any $g \in C(G(n,m))$ we may write:
\[
R_{m}^*(g)(\theta) = \int_{G_\theta(n,m)} g(E) d\nu_\sigma(E).
\]

We will need the following fact, which is an immediate corollary
of Proposition \ref{prop:Grassmann}. We postpone the formulation
and proof of Proposition \ref{prop:Grassmann} for the Appendix,
as the technique involved is different in spirit to the rest of
this note.

\begin{cor} \label{cor:Grassmann}
Let $n>1$ and let $k_1,k_2 \geq 1$ denote integers such that $l =
k_1 + k_2 \leq n-1$. Let $\theta \in S^{n-1}$. For $a=k_1,k_2,l$,
denote by $G^{a} = G(n,n-a)$ and by $\mu^{a}_\theta$ the Haar
probability measure on $G^a_\theta$. For $F \in G^l$ and
$a=k_1,k_2$, denote by $\mu^{a}_F$ the Haar probability measure
on $G^a_F$. Then for any continuous function $f(E_1,E_2)$ on
$G^{k_1}\times G^{k_2}$:
%Denote by $\bar{E} =
%(E_1,\ldots,E_r)$ an ordered set with $E_i \in G^{k_i}$.
\begin{eqnarray}
\nonumber & & \!\!\!\!\!\!\!\!\! \int_{E_1 \in G^{k_1}_\theta}
\int_{E_2\in G^{k_2}_\theta}
f(E_1,E_2) d\mu^{k_1}_\theta(E_1) d\mu^{k_2}_\theta(E_2) = \\
\nonumber & & \!\!\!\!\!\!\!\!\! \int_{F \in G^l_\theta} \int_{E_1
\in G^{k_1}_F}\int_{E_2 \in G^{k_2}_F} f(E_1,E_2) \Delta(E_1,E_2)
d\mu^{k_1}_F(E_1) d\mu^{k_2}_F(E_2) d\mu^{l}_\theta(F),
\end{eqnarray}
where $\Delta(E_1,E_2)$ is some (known) non-negative continuous
function on $G^{k_1} \times G^{k_2}$.
\end{cor}

We will show the following basic property of $k$-Busemann-Petty
bodies, and immediately deduce (3a), (3b) and (3c) from the
Structure Theorem in the Introduction.

\begin{prop} \label{prop:BP-basic-prop}
Let $K_1 \in \BP_{k_1}^n$ and $K_2 \in \BP_{k_2}^n$ for $k_1,k_2
\geq 1$ such that $l = k_1 + k_2 \leq n-1$. Then the star-body $L$
defined by $\rho_L^{l} = \rho_{K_1}^{k_1} \rho_{K_2}^{k_2}$
satisfies $L \in \BP_{l}^n$.
\end{prop}

\begin{proof}
First, assume that $K_i \in \BP_{k_i}^{n,\infty}$ for $i=1,2$, so
that $\rho_K^{k_i} = R_{n-k_i}^*(g_i)$ with $g_i \in
C_+^\infty(G(n,n-k_i))$. Using the notations and result of
Corollary \ref{cor:Grassmann}, we have:
\begin{eqnarray}
\nonumber \rho_L^{l}(\theta) = \rho_{K_1}^{k_1}(\theta)
\rho_{K_2}^{k_2}(\theta) = \int_{E_1 \in G^{k_1}_\theta} g_1(E_1)
d\mu^{k_1}_\theta(E_1) \int_{E_2 \in G^{k_2}_\theta} g_2(E_2)
d\mu^{k_2}_\theta(E_2) \\
\nonumber = \int_{F \in G^l_\theta} \int_{E_1 \in
G^{k_1}_F}\int_{E_2 \in G^{k_2}_F} g(E_1)g(E_2) \Delta(E_1,E_2)
d\mu^{k_1}_F(E_1) d\mu^{k_2}_F(E_2) d\mu^{l}_\theta(F).
\end{eqnarray}
Denoting:
\[
h(F) = \int_{E_1 \in G^{k_1}_F}\int_{E_2 \in G^{k_2}_F}
g(E_1)g(E_2) \Delta(E_1,E_2) d\mu^{k_1}_F(E_1) d\mu^{k_2}_F(E_2),
\]
we see that $h(F)$ is a non-negative continuous function on
$G(n,n-l)$. Therefore:
\[
\rho_L^{l}(\theta) = \int_{F \in G^l_\theta} h(F)
d\mu^{l}_\theta(F),
\]
implying that $L \in \BP_l^n$. The general case, when $K_i \in
\BP_{k_i}^n$ without any smoothness assumptions, follows from
Remark \ref{rem:BP-use-smooth}. Indeed, by approximating each
$K_i$ in the radial metric by smooth bodies $\set{K^m_i} \subset
\BP_{k_i}^n$, the bodies $\set{L^m}$ defined by $\rho_{L^m}^{l} =
\rho_{K^m_1}^{k_1} \rho_{K^m_2}^{k_2}$ satisfy that $L^m \in
\BP_{l}^n$ and obviously $L^m$ approximate $L$ in the radial
metric, implying that $L \in \BP_l^n$.
\end{proof}

Applying Proposition \ref{prop:BP-basic-prop} with $K_1 = K_2$,
we have:
\begin{cor} \label{cor:BP-k1-and-k2}
$\BP_{k_1}^n \cap \BP_{k_2}^n \subset \BP_{k_1+k_2}^n$ for
$k_1,k_2 \geq 1$ such that $k_1 + k_2 \leq n-1$.
\end{cor}

By successively applying Corollary \ref{cor:BP-k1-and-k2}, we see
that $\BP_k^n \subset \BP_l^n$ if $k$ divides $l$. The question
whether $\BP_k^n \subset \BP_l^n$ for general $1\leq k < l \leq
n-1$ remains open. Nevertheless, we are able to show the
following "non-linear" embedding of $\BP_k^n$ into $\BP_l^n$,
which is again an immediate corollary of Proposition
\ref{prop:BP-basic-prop} (using $K_2 = D_n \in \BP_{l-k}^n$):

\begin{prop} \label{prop:BP-non-linear-embedding}
If $K \in \BP_{k}^n$ then the star-body $L$ defined by $\rho_L =
\rho_K^{k/l}$ satisfies $L \in \BP_{l}^n$ for $1\leq k \leq l \leq
n-1$.
\end{prop}

We prefer to give another proof of this statement, one which does
not rely on Proposition \ref{prop:Grassmann}.

\begin{proof}
Assume that $K \in \BP_k^{n,\infty}$, so that $\rho_K^k =
R_{n-k}^*(g_K)$ and $g_K \in C_+^\infty(G(n,n-k))$, and define the
star-body $L$ by $\rho_L = \rho_K^{k/l}$. For $\theta \in
S^{n-1}$ and $a = k,l$ denote by $\mu^a_\theta$ the Haar
probability measure on $G_\theta(n,n-a)$. For $F \in G(n,n-l)$,
denote by $\mu^k_F$ the Haar probability measure on $G_F(n,n-k)$.
Then:
\begin{eqnarray}
\nonumber \rho_L^l(\theta) = \rho_K^k(\theta) = \int_{G_\theta(n,n-k)} g_K(E) d\mu^k_{\theta}(E) = \\
\nonumber \int_{G_\theta(n,n-l)} \int_{G_F(n,n-k)} g_K(E)
d\mu^k_F(E) d\mu^l_{\theta}(F).
\end{eqnarray}
The last transition is justified by the fact that the probability
measure $d\mu^k_F(E) d\mu^l_{\theta}(F)$ on $G_\theta(n,n-k)$ is
invariant under orthogonal rotations in $O(n)$ which preserve
$\theta$, and therefore coincides with $d\mu^k_{\theta}(E)$, the
Haar probability measure on $G_\theta(n,n-k)$. Defining $g_L \in
C_+(G(n,n-l))$ by $g_L(F) = \int_{G_F(n,n-k)} g(E) d\mu^k_F(E)$
for $F \in G(n,n-l)$, we see that:
\[
\rho_L^l(\theta) = R_{n-l}^*\brac{g_L}(\theta).
\]
Together with Remark \ref{rem:BP-use-smooth}, this concludes the
proof.
\end{proof}

The Ellipsoid Corollary from the Introduction should now be
clear. We repeat it here for convenience:
\begin{cor}
For any $1 \leq k \leq n-1$ and $k$ ellipsoids
$\set{\E_i}_{i=1}^k$ in $\Real^n$, define the body $L$ by:
\[
\rho_L =  \rho_{\E_1} \cdot \ldots \cdot \rho_{\E_k},
\]
and let $k\leq l\leq n-1$. Then there exists a sequence of
star-bodies $\set{L_i}$ which tends to $L$ in the radial metric
and satisfies:
\[
\rho_{L_i} = \rho_{\E^i_1}^{l} + \ldots + \rho_{\E^i_{m_i}}^{l},
\]
where $\set{\E^i_j}$ are ellipsoids.
\end{cor}
\begin{proof}
The body $L_2$ defined by $\rho^k_{L_2} = \rho_L$ is in $\BP_k^n$
by Proposition \ref{prop:BP-basic-prop} (applied successively to
the ellipsoids $\set{\E_i}$, which are in $\BP_1^n$). For $l>k$,
Proposition \ref{prop:BP-non-linear-embedding} implies that the
body $L_3$ defined by $\rho^l_{L_3} = \rho^k_{L_2} = \rho_L$ is in
$\BP_l^n$, otherwise this is trivial. Using Grinberg and Zhang's
characterization of $\BP_l^n$ (Theorem \ref{thm:G&Z}), the claim
is established.
\end{proof}

Incidentally, Proposition \ref{prop:BP-basic-prop} also enables us
to give the following strengthened version of Theorem
\ref{thm:G&Z}:
\begin{cor}
A star-body $K$ is a $k$-Busemann-Petty body iff it is the limit
of $\set{K_i}$ in the radial metric, where each $K_i$ is of the
following form:
\[
\rho^k_{K_i} = \rho_{\E^{i}_{1,1}} \cdot \ldots \cdot
\rho_{\E^{i}_{1,k}} + \ldots + \rho_{\E^{i}_{m_i,1}} \cdot \ldots
\cdot \rho_{\E^{i}_{m_i,k}},
\]
where $\set{\E^{i}_{j,l}}$ are ellipsoids.
\end{cor}
\begin{proof}
Obviously this representation generalizes the one given by
Grinberg and Zhang in Theorem \ref{thm:G&Z}, so it is enough to
show the "if" part. But this follows from the closure of
$\BP_k^n$ under limit in the radial metric, $k$-radial sums, and
Proposition \ref{prop:BP-basic-prop} (which as above shows that
the body $L$ defined by $\rho_L^k = \rho_{\E_1} \cdot \ldots
\cdot \rho_{\E_k}$ is in $\BP_k^n$).
\end{proof}

For completeness, we conclude our investigation of the structure
of $\BP_k^n$ with the following result of Grinberg and Zhang from
\cite{Grinberg-Zhang}. Their argument is the same one used by
Goodey and Weil for intersection bodies ($\BP_1^n$), and is an
immediate corollary of Theorem \ref{thm:G&Z}.
\begin{cor}[Grinberg and Zhang]
If $K \in \BP_k^n$ then any $m$-dimensional central section $L$ of
$K$ (for $m>k$) satisfies $L \in \BP_k^m$.
\end{cor}
\begin{proof}
Since and central section of an ellipsoid is again an ellipsoid,
the claim follows immediately from Theorem \ref{thm:G&Z}.
\end{proof}

%%%%%%%%%%%%%%%%%%%%%%%%%%%%%%%%%%%%%%%%%%%%%%%%%%%%%%%%%%%%%%%%%%%%%%%%%%%%%%%%%%%%%%%%%%%%%
%%%%%%%%%%%%%%%%%%%%%%%%%%%%%%%%%%%%%%%%%%%%%%%%%%%%%%%%%%%%%%%%%%%%%%%%%%%%%%%%%%%%%%%%%%%%%
%%%%%%%%%%%%%%%%%%%%%%%%%%%%%%%%%%%%%%%%%%%%%%%%%%%%%%%%%%%%%%%%%%%%%%%%%%%%%%%%%%%%%%%%%%%%%
%%%%%%%%%%%%%%%%%%%%%%%%%%%%%%%%%%%%%%%%%%%%%%%%%%%%%%%%%%%%%%%%%%%%%%%%%%%%%%%%%%%%%%%%%%%%%

\bigskip

We now turn to prove the Structure Theorem from the Introduction
for $\I_k^n$. As will be evident, the techniques involved are
totally different from those which were used for $\BP_k^n$. The
only point of similarity is Lemma \ref{lem:I-dense} below. We
denote by $\I_k^{n,\infty}$ the class of infinitely smooth
$k$-intersection bodies in $\Real^n$. As mentioned in Section
\ref{sec:1}, this implies for $K \in \I_k^{n,\infty}$ that
$\norm{\cdot}_K^{-k} , (\norm{\cdot}_K^{-k})^\wedge \in
C^\infty(\Real^n \setminus \set{0})$. We begin with the following
useful lemma:

\begin{lem} \label{lem:ast-commutes-with-wedge}
For any $p \in (-n,0)$, $g \in C^\infty(S^{n-1})$ and $\mu \in
\M(O(n))$, $E_p^\wedge(\mu \ast g) = \mu \ast E_p^\wedge(g)$ as
functions on $\Real^n \setminus \set{0}$.
\end{lem}
\begin{proof}
First, let us extend the definition of $\mu \ast f$ to any
function $f \in C(\Real^n)$, as follows: $(\mu \ast f)(x) =
\int_{O(n)} f(u(x)) d\mu(u)$ for every $x\in\Real^n$. Next,
notice that for a test function $\phi$, $(\mu \ast \phi)^\wedge =
\mu \ast \hat{\phi}$. Indeed, when $\mu$ is a delta function at
$u \in O(n)$, $(\phi(u(\cdot)))^\wedge(x) = \hat{\phi}(u(x))$
because the Fourier transform commutes with the action of $O(n)$.
And for a general $\mu \in \M(O(n))$, by Fubini's Theorem:
\begin{eqnarray}
\nonumber (\mu \ast \phi)^\wedge(x) &= & \int_{\Real^n}
\int_{O(n)}
\phi(u(y)) d\mu(u) \exp(-i\scalar{y,x}) dy \\
\nonumber &= & \int_{O(n)}
\int_{\Real^n} \phi(u(y)) \exp(-i\scalar{y,x}) dy d\mu(u) \\
\nonumber & = & \int_{O(n)} (\phi(u(\cdot)))^\wedge(x) d\mu(u) =
\int_{O(n)} \hat{\phi}(u(x)) d\mu(u) = \mu \ast \hat{\phi}.
\end{eqnarray}
Since $g, \mu \ast g \in C^\infty(S^{n-1})$, it follows that
$E_p^\wedge(\mu \ast g), \mu \ast E_p^\wedge(g) \in
C^\infty(\Real^n \setminus \set{0})$, and for any test function
$\phi$:
\begin{eqnarray}
\nonumber \scalar{E_p^\wedge(\mu \ast g),\phi} = \scalar{E_p(\mu
\ast g),\hat{\phi}} = \scalar{\mu \ast E_p(g),\hat{\phi}} =
\scalar{E_p(g),\mu^{-1} \ast \hat{\phi}} \\
\nonumber = \scalar{E_p(g),(\mu^{-1} \ast \phi)^\wedge} =
\scalar{E_p^\wedge(g),\mu^{-1} \ast \phi} = \scalar{\mu \ast
E_p^\wedge(g),\phi}.
\end{eqnarray}
Therefore $E_p^\wedge(\mu \ast g) = \mu \ast E_p^\wedge(g)$ as
functions.
\end{proof}

\begin{lem} \label{lem:I-dense}
$\I_k^{n,\infty}$ is dense in $\I_k^n$.
\end{lem}
\begin{proof}
Let $K \in \I_k^n$, and let $\mu \in \M_+(S^{n-1})$ be the
measure from Proposition \ref{prop:I-spherical-char} satisfying
for every $f\in C^\infty(S^{n-1})$:
\[
\int_{S^{n-1}} f(\theta) \rho_K^k(\theta) d\sigma(\theta) =
\int_{S^{n-1}} E_{-n+k}^\wedge(f)(\theta) d\mu(\theta).
\]
Let $\set{u_i} \subset C^\infty(O(n))$ be an approximate identity
as in Lemma \ref{lem:approximate-identity}, and let $K_i$ be the
star-body for which $\rho_{K_i}^k = u_i \ast \rho_{K}^k$. Then by
Lemma \ref{lem:approximate-identity}, $\set{K_i}$ is a sequence of
infinitely smooth star-bodies which tend to $K$ in the radial
metric. It remains to check that each $K_i$ is a $k$-intersection
body. Indeed, using the notations of Section \ref{sec:1} and Lemma
\ref{lem:ast-commutes-with-wedge}, for any $f\in
C^\infty(S^{n-1})$:
\begin{eqnarray}
\nonumber \scalar{f,\rho_{K_i}^k} & = & \scalar{f,u_i \ast
\rho_{K}^k} = \scalar{u_i^{-1} \ast f,\rho_K^k} =
\scalar{E_{-n+k}^\wedge(u_i^{-1} \ast f),\mu} \\
\nonumber & = & \scalar{u_i^{-1} \ast E_{-n+k}^\wedge(f),\mu} =
\scalar{E_{-n+k}^\wedge(f),u_i \ast \mu}.
\end{eqnarray}
Since $u_i \ast \mu \in C^\infty_+(S^{n-1})$, again by
Proposition \ref{prop:I-spherical-char} this implies that $K_i
\in \I_k^n$.
\end{proof}

\begin{rem} \label{rem:I-use-smooth}
By the Lemma and the closure of $\I_k^n$ (for any $k =
1,\ldots,n-1$) under limit in the radial metric, it is enough to
prove all the remaining items for the classes $\I_k^{n,\infty}$.
\end{rem}

For the next fundamental proposition, we will need the following
observation. It is classical that for two test functions
$\phi_1,\phi_2$, $(\phi_1 \phi_2)^\wedge = \hat{\phi_1} \star
\hat{\phi_1}$ where $\star$ denotes the standard convolution on
$\Real^n$. In general, the convolution of two distributions does
not exist. Nevertheless, when the two distributions $f_1$,$f_2$
are locally integrable homogeneous functions with the right
degrees, their convolution may be defined as usual. Assume that
$f_i$ is even homogeneous of degree $-n+p_i$ for $p_i > 0$ and
that $p_1 + p_2 < n$. Since $f_i$ are locally integrable and at
infinity their product decays faster than $\abs{x}^{-n}$, the
following integral converges for $x \in \Real^n \setminus
\set{0}$:
\begin{equation} \label{eq:conv-defn}
f_1 \star f_2 (x) = \int f_1(x-y) f_2(y) dy.
\end{equation}
It is easy to check that with this definition, $f_1 \star f_2$ is
homogeneous of degree $-n+p_1+p_2$, hence again locally
integrable. Now assume in addition that $f_i$ are infinitely
smooth functions on $\Real^n \setminus \set{0}$, and therefore so
are $\hat{f_i}$. We claim that as distributions $(f_1 \star
f_2)^\wedge = \hat{f_1} \hat{f_2}$. To see this, we define the
product and convolution of an even distribution $f$ with an even
test-function $\phi$, as the distributions denoted $\phi f$ and
$\phi \star f$, respectively, satisfying for any test function
$\varphi$ that:
\[
\nonumber \scalar{\phi f , \varphi} = \scalar{f , \phi \varphi}
\text{ and }\scalar{\phi \star f , \varphi} = \scalar{f , \phi
\star \varphi}.
\]
When $f$ is a locally integrable function, it is clear that $\phi
f$ and $\phi \star f$ as distributions coincide with the usual
product and convolution as functions. The same reasoning shows
that when $f_1,f_2$ are locally integrable even functions such
that $f_1f_2$ is integrable at infinity (as before the definition
in (\ref{eq:conv-defn})), we have:
\begin{equation} \label{eq:conv1}
\scalar{f_1 \star f_2, \phi} = \scalar{f_1, \phi \star f_2},
\end{equation}
where the action $\scalar{\cdot,\cdot}$ is interpreted here and
henceforth as integration in $\Real^n$. Similarly, when $f_1 f_2$
is locally integrable, we have:
\begin{equation} \label{eq:conv1b}
\scalar{f_1 f_2 , \phi} = \scalar{f_1, \phi f_2}.
\end{equation}

With the above definitions, we see that $(\phi \star f)^\wedge =
\hat{\phi} \hat{f}$ because for any test function $\varphi$:
\begin{equation} \label{eq:conv2}
\scalar{\phi \star f , \hat{\varphi}} = \scalar{f , \phi \star
\hat{\varphi}} = \scalar{\hat{f} , \hat{\phi} \varphi} =
\scalar{\hat{\phi} \hat{f} , \varphi}.
\end{equation}
Now when $f,g$ are two locally integrable infinitely smooth
functions on $\Real^n \setminus \set{0}$, such that $\hat{f} g$
is locally integrable, it is easy to see that we may replace
$\varphi$ in (\ref{eq:conv2}) with $g$. The reason is that we may
weakly approximate $g$ with test functions $g_i$ such that $\int
h g_i \rightarrow \int h g$ and $\int h \hat{g_i} \rightarrow
\int h \hat{g}$, for any locally integrable continuous function
$h$ on $\Real^n \setminus \set{0}$ such that $\int h g$ exists.
For instance, we may use $g_i = (g \star \delta_i)
\hat{\delta_i}$, where $\delta_i$ are Gaussians tending to a
delta-function at 0; by (\ref{eq:conv2}) it is clear that
$\hat{g_i} = ( \hat{g} \hat{\delta_i} ) \star \delta_i$, which
weakly tends to $\hat{g}$ (by testing against a test-function). We
summarize this by writing:
\begin{equation} \label{eq:conv3}
\scalar{\phi \star f , \hat{g}} = \scalar{\hat{\phi} \hat{f}, g}.
\end{equation}

Combining (\ref{eq:conv1}), (\ref{eq:conv1b}) and (\ref{eq:conv3})
and using the fact that $f_i,\hat{f_i},\hat{f_1}\hat{f_2}$ are
infinitely smooth and locally integrable, we see that for any
even test function $\phi$:
\[
\scalar{(f_1 \star f_2)^\wedge,\phi} = \scalar{f_1 \star f_2,
\hat{\phi}} = \scalar{f_1, \hat{\phi} \star f_2} =
%1/(2\pi)^n \scalar{\hat{f_1},(\hat{\phi} \star f_2)^\wedge} =
\scalar{\hat{f_1},\phi \hat{f_2}} = \scalar{\hat{f_1} \hat{f_2} ,
\phi}.
\]
This proves that under the above conditions:
\begin{equation} \label{eq:convolution-identity}
(f_1 \star f_2)^\wedge = \hat{f_1} \hat{f_2}.
\end{equation}

\begin{rem}
Note that the homogeneity of $f_1,f_2$ was not used, we only
needed the appropriate asymptotic behaviour at 0 and infinity.
Using the homogeneity, a different approach to derive
(\ref{eq:convolution-identity}) was suggested to us by A.
Koldobsky, by applying \cite[Lemma 1]{Koldobsky-convolution}.
With this approach, the smoothness assumptions on $f_1,f_2$ may
be omitted, and (\ref{eq:convolution-identity}) is understood as
equality between distributions.
\end{rem}

\medskip

Using this notion of convolution, we can now show the following
basic property of $k$-intersection bodies, and immediately deduce
(3a), (3b) and (3c) from the Structure Theorem in the
Introduction. The following was also recently noticed
independently by Koldobsky (but not published):

\begin{prop} \label{prop:I-basic-prop}
Let $K_1 \in \I_{k_1}^n$ and $K_2 \in \I_{k_2}^n$ for $k_1,k_2
\geq 1$ such that $l = k_1 + k_2 \leq n-1$. Then the star-body $L$
defined by $\rho_L^{l} = \rho_{K_1}^{k_1} \rho_{K_2}^{k_2}$
satisfies $L \in \I_{l}^n$.
\end{prop}

\begin{proof}
First, assume that $K_i \in \I_{k_i}^{n,\infty}$ for $i=1,2$, so
that $(\norm{\cdot}_K^{-k_i})^\wedge \in C^\infty_+(\Real^n
\setminus \set{0})$ and is homogeneous of degree $-n+k_i$. Since
$l < n$ the convolution $(\norm{\cdot}_K^{-k_1})^\wedge \star
(\norm{\cdot}_K^{-k_2})^\wedge$ as distributions is well defined
(as explained above). Therefore:
\[
(\norm{\cdot}_L^{-l})^\wedge = (\norm{\cdot}_K^{-k_1}
\norm{\cdot}_K^{-k_2})^\wedge = (\norm{\cdot}_K^{-k_1})^\wedge
\star (\norm{\cdot}_K^{-k_2})^\wedge \geq 0,
\]
as a function on $\Real^n \setminus \set{0}$, which implies that
$L \in \I_l^n$. The general case, when $K_i \in \I_{k_i}^n$
without any smoothness assumptions, follows from Remark
\ref{rem:I-use-smooth} in the same manner as in the proof of
Proposition \ref{prop:BP-basic-prop}.
\end{proof}

Applying Proposition \ref{prop:I-basic-prop} with $K_1 = K_2$, we
have:
\begin{cor} \label{cor:I-k1-and-k2}
$\I_{k_1}^n \cap \I_{k_2}^n \subset \I_{k_1+k_2}^n$ for $k_1,k_2
\geq 1$ such that $k_1 + k_2 \leq n-1$.
\end{cor}

By successively applying Corollary \ref{cor:I-k1-and-k2}, we see
that $\I_k^n \subset \I_l^n$ if $k$ divides $l$. As for the class
$\BP$, the question whether $\I_k^n \subset \I_l^n$ for general
$1\leq k < l \leq n-1$ remains open. Nevertheless, we are able to
show again the following "non-linear" embedding of $\I_k^n$ into
$\I_l^n$, which is again an immediate corollary of Proposition
\ref{prop:I-basic-prop} (using $K_2 = D_n \in \I_{l-k}^n$):

\begin{cor} \label{cor:I-non-linear-embedding}
If $K \in \I_{k}^n$ then the star-body $L$ defined by $\rho_L =
\rho_K^{k/l}$ satisfies $L \in \I_{l}^n$ for $1\leq k \leq l \leq
n-1$.
\end{cor}

We conclude this section with our last observation:

\begin{prop}
If $K \in \I_k^n$ then any $m$-dimensional central section $L$ of
$K$ (for $m>k$) satisfies $L \in \I_k^m$.
\end{prop}

\begin{proof}
Let $K$ be a star-body in $\Real^n$, fix $k \in
\set{1,\ldots,n-2}$, and let $H \in G(n,m)$ for $m>k$. In view of
Theorem \ref{thm:I-char}, we have to show that as distributions:
\[
(\norm{\cdot}_K^{-k})^\wedge \geq 0 \; \textnormal{ implies } \;
(\norm{\cdot}_K^{-k} |_H)^\wedge = (\norm{\cdot}_{K \cap
H}^{-k})^\wedge \geq 0.
\]
This becomes intuitively clear, after noticing that for a test
function $\phi$:
\[
(\phi|_H)^\wedge (u) = \int_{u + H^\perp} \hat{\phi}(y) dy.
\]
Nevertheless, for a more general function $f =
\norm{\cdot}_K^{-k}$ such that $\hat{f} \geq 0$ as a distribution,
we will need a somewhat different proof. Note that since $m>k$,
$f$ is locally integrable on any affine translate $z + H$, and
that for any test function $\phi_H$ on $H$, $\int_{H} f(y+z)
\phi(y) dy$ is continuous w.r.t. $z \in H^\perp$. Now let
$\phi_H$ be any \emph{non-negative} test function on $H$. For
$\epsilon>0$, denote by $\varphi_{H^\perp,\epsilon}$ the
(positive) Gaussian function on $H^\perp$ such that
$(\varphi_{H^\perp,\epsilon})^\wedge$ is the density function of
a standard Gaussian variable on $H^\perp$ with covariance matrix
$\epsilon I_{H^\perp}$. For $y \in H$ and $z \in H^\perp$, define
$\phi_\epsilon(y+z) = \phi_H(y) \varphi_{H^\perp,\epsilon}(z)$.
Clearly $\phi_\epsilon$ is a test function on $\Real^n$,
$\phi_\epsilon \geq 0$, and $(\phi_\epsilon)^\wedge(y+z) =
(\phi_H)^\wedge(y) (\varphi_{H^\perp,\epsilon})^\wedge(z)$. We
therefore have:
\begin{eqnarray}
\nonumber \scalar{(f|_H)^\wedge,\phi_H} & = &
\scalar{f|_H,(\phi_H)^\wedge} = \int_H f(y) (\phi_H)^\wedge(y) dy
\\
\nonumber & = & \lim_{\epsilon \rightarrow 0} \int_{H^\perp}
(\varphi_{H^\perp,\epsilon})^\wedge(z) \int_H
f(y+z) (\phi_H)^\wedge(y) dy \; dz \\
\nonumber & = & \lim_{\epsilon \rightarrow 0} \int_{\Real^n} f(x)
(\phi_\epsilon)^\wedge(x) dx  \\
\nonumber & = & \lim_{\epsilon \rightarrow 0}
\scalar{f,(\phi_\epsilon)^\wedge}
 = \lim_{\epsilon \rightarrow 0} \scalar{\hat{f},\phi_\epsilon} \geq 0.
\end{eqnarray}
Since $\phi_H \geq 0$ was arbitrary, it follows that
$(f|_H)^\wedge \geq 0$.
\end{proof}

%%%%%%%%%%%%%%%%%%%%%%%%%%%%%%%%%%%%%%%%%%%%%%%%%%%%%%%%%%%%%%%%%%%%%%%%%%%%%%%%%%%%%%%
%%%%%%%%%%%%%%%%%%%%%%%%%%%%%%%%%%%%%%%%%%%%%%%%%%%%%%%%%%%%%%%%%%%%%%%%%%%%%%%%%%%%%%%
%%%%%%%%%%%%%%%%%%%%%%%%%%%%%%%%%%%%%%%%%%%%%%%%%%%%%%%%%%%%%%%%%%%%%%%%%%%%%%%%%%%%%%%
%%%%%%%%%%%%%%%%%%%%%%%%%%%%%%%%%%%%%%%%%%%%%%%%%%%%%%%%%%%%%%%%%%%%%%%%%%%%%%%%%%%%%%%
%%%%%%%%%%%%%%%%%%%%%%%%%%%%%%%%%%%%%%%%%%%%%%%%%%%%%%%%%%%%%%%%%%%%%%%%%%%%%%%%%%%%%%%
%%%%%%%%%%%%%%%%%%%%%%%%%%%%%%%%%%%%%%%%%%%%%%%%%%%%%%%%%%%%%%%%%%%%%%%%%%%%%%%%%%%%%%%
%%%%%%%%%%%%%%%%%%%%%%%%%%%%%%%%%%%%%%%%%%%%%%%%%%%%%%%%%%%%%%%%%%%%%%%%%%%%%%%%%%%%%%%
%%%%%%%%%%%%%%%%%%%%%%%%%%%%%%%%%%%%%%%%%%%%%%%%%%%%%%%%%%%%%%%%%%%%%%%%%%%%%%%%%%%%%%%

\section{The connection between Radon and Fourier Transforms}
\label{sec:3}

We have seen that although the classes $\BP_k^n$ and $\I_k^n$
share the exact same structure and easily verify that $\BP_k^n
\subset \I_k^n$, they are defined and handled using very
different notions: Radon and Fourier transforms, respectively.
The aim of this section is to establish a common ground that will
enable to attack the question of whether $\BP_k^n = \I_k^n$ from a
unified point of view. Since $\BP_k^n \subset \I_k^n$, it seems
natural that this common ground will involve the language of
Radon transforms, so we will have to translate the action of the
Fourier transform to this language.

We will use the following notation. If $\mu \in \M(G(n,n-m))$, we
denote by $\mu^\perp \in \M(G(n,m))$ the measure defined by
$\mu^\perp(A) = \mu(A^\perp)$ for any Borel set $A \subset
G(n,m)$, where $A^\perp = \set{E^\perp | E \in A}$. Note that the
operation $\mu \rightarrow \mu^\perp$ is dual to the operator
$I:C(G(n,m))\rightarrow C(G(n,n-m))$ defined in the Introduction,
in the sense that $\scalar{\mu,I(f)}_{G(n,n-m)} =
\scalar{\mu^\perp,f}_{G(n,m)}$. We recall that $I(f)(E) =
f^\perp(E) = f(E^\perp)$ for any $E\in G(n,n-m)$. We therefore
extend $I$ to an operator $I:\M(G(n,m)) \rightarrow
\M(G(n,n-m))$, defined as $I(\mu) = \mu^\perp$, and by abuse of
notation we say that $I$ is self-dual.

Theorem \ref{thm:R_k(wedge)} in Section \ref{sec:1} was the first
example relating the Radon and Fourier transforms. Using operator
notations, this may be stated as:
\begin{equation} \label{eq:R_k(wedge)-operator}
R_{n-k} \circ E_{-k}^\wedge =  c(n,k) I \circ R_k \; ,
\end{equation}
as operators from $C^\infty_e(S^{n-1})$ to
$C^\infty_e(G(n,n-k))$. In view of the remark immediately after
Theorem \ref{thm:R_k(wedge)}, a generalization of
(\ref{eq:R_k(wedge)-operator}) is given by the Spherical Parseval
identity, which in the formulation of Remark
\ref{rem:spherical-parseval}, shows that
%A more general property
%of the Fourier transform was given by the Spherical Parseval
%identity, which in the formulation of Remark \ref{rem:spherical-parseval}, show that
$E_{-k}^\wedge$ is a self-adjoint operator on
$C_e^\infty(S^{n-1})$:
\begin{equation} \label{eq:SphericalParseval-operator}
(E_{-k}^\wedge)^* = E_{-k}^\wedge.
\end{equation}
Passing to the dual in (\ref{eq:R_k(wedge)-operator}) and using
(\ref{eq:SphericalParseval-operator}), we immediately have:
\begin{equation} \label{eq:wedge(R_k^*)}
E_{-k}^\wedge \circ R_{n-k}^* = c(n,k) R_k^* \circ I \; ,
\end{equation}
as operators on certain spaces. We formulate this more carefully
in the next Proposition:

\begin{prop} \label{prop:wedge(R_k^*)}
Let $f \in C^\infty_e(S^{n-1})$, and assume that $f =
R_{n-m}^*(d\mu)$ as measures in $\M(S^{n-1})$, for some measure
$\mu \in \M(G(n,n-m))$. Then $E_{-m}^\wedge(f) = c(n,m)
R_m^*(d\mu^\perp)$ as measures in $\M(S^{n-1})$, where $c(n,m)$ is
the constant from Lemma \ref{lem:E_p(1)}.
\end{prop}
\begin{proof}
Let $g \in C^\infty_e(S^{n-1})$ be arbitrary. Then by the
Spherical Parseval identity and Theorem \ref{thm:R_k(wedge)}:
\begin{eqnarray}
\nonumber & & \int_{S^{n-1}} E_{-m}^\wedge(f)(\theta) g(\theta)
d\sigma(\theta) = \int_{S^{n-1}} f(\theta)
E_{-m}^\wedge(g)(\theta) d\sigma(\theta) \\
\nonumber & = & \int_{S^{n-1}} R_{n-m}^*(d\mu)(\theta)
E_{-m}^\wedge(g)(\theta) d\sigma(\theta) = \int_{G(n,n-m)}
R_{n-m}(E_{-m}^\wedge(g))(F) d\mu(F) \\
\nonumber & = & c(n,m) \int_{G(n,n-m)} R_{m}(g)(F^\perp) d\mu(F) =
c(n,m) \int_{G(n,m)} R_{m}(g)(F) d\mu^\perp(F) \\
\nonumber & = & c(n,m) \int_{S^{n-1}}
R_{m}^*(d\mu^\perp)(\theta)g(\theta)d\sigma(\theta).
\end{eqnarray}
Since $C^\infty_e(S^{n-1})$ is dense in $C_e(S^{n-1})$ in the
maximum norm, the proposition follows.
\end{proof}

\medskip

In the context of star-bodies, the following is an immediate
corollary of Proposition \ref{prop:wedge(R_k^*)}:
\begin{cor} \label{cor:wedge(R_k^*)}
Let $K$ be an infinitely smooth star-body in $\Real^n$. Then for
a measure $\mu \in \M(G(n,n-k))$:
\[
\norm{\cdot}_K^{-k} = R_{n-k}^*(d\mu) \text{  iff   }
(\norm{\cdot}_K^{-k})^\wedge = c(n,k) R_{k}^*(d\mu^\perp),
\]
where $c(n,k)$ is the constant from Lemma \ref{lem:E_p(1)}, and
the equalities are understood as equalities between measures in
$\M(S^{n-1})$.
\end{cor}
\begin{proof}
The "only if" part follows immediately from Proposition
\ref{prop:wedge(R_k^*)} with $m=k$ and $f = \norm{\cdot}_K^{-k}$.
The "if" part follows by applying Proposition
\ref{prop:wedge(R_k^*)} with $m=n-k$ and $f =
(\norm{\cdot}_K^{-k})^\wedge$, and using the fact that
$E_{-n+k}^\wedge(f) = (2\pi)^n \norm{\cdot}_K^{-k}$ and that the
constants $c(n,k)$ from Lemma \ref{lem:E_p(1)} satisfy $c(n,k)
c(n,n-k) = (2\pi)^n$.
\end{proof}

Proposition \ref{prop:wedge(R_k^*)} has several interesting
consequences. The first one is:

\begin{thm}
Let $n>1$ and fix $1\leq k \leq n-1$. Then:
\[\BP_k^n = \I_k^n \text{  iff  } \; \BP_{n-k}^n = \I_{n-k}^n.
\]
\end{thm}
\begin{proof}
Assume that $\BP_{n-k}^n = \I_{n-k}^n$, and let $K \in
\I_{n-k}^{n,\infty}$. In view of Lemma \ref{lem:I-dense}, the fact
that $\BP_k^n$ is closed under limit in the radial metric, and
Corollary \ref{cor:BP-subset-I}, it is enough to show that $K \in
\BP_k^n$. Since $(\norm{\cdot}_K^{-k})^\wedge \geq 0$ by Theorem
\ref{thm:I-char}, we may define the infinitely smooth star-body
$L$ as the body for which $\norm{\cdot}_L^{-n+k} =
(\norm{\cdot}_K^{-k})^\wedge$. Therefore
$(\norm{\cdot}_L^{-n+k})^\wedge = (2\pi)^n \norm{\cdot}_K^{-k}
\geq 0$, hence $L \in \I_{n-k}^n$. It follows from our assumption
that $L \in \BP_{n-k}^n$, so there exists a non-negative measure
$\mu \in \M_+(G(n,k))$ so that $(\norm{\cdot}_K^{-k})^\wedge =
\norm{\cdot}_L^{-n+k} = R_k^*(d\mu)$. By Corollary
\ref{cor:wedge(R_k^*)}, this implies that $\norm{\cdot}_K^{-k} =
c(n,k) R_{n-k}^*(d\mu^\perp)$. Therefore $K \in \BP_k^n$, which
concludes the proof.
\end{proof}

Another immediate consequence of Proposition
\ref{prop:wedge(R_k^*)} is another elementary proof of:
\begin{cor} \label{cor:BP-subset-I-again}
\[
\BP_k^n \subset \I_k^n.
\]
\end{cor}
\begin{proof}
Let $K \in \BP_k^{n,\infty}$, so $\norm{\cdot}_K^{-k} =
R_{n-k}^*(d\mu)$ for some non-negative Borel measure $\mu \in
\M_+(G(n,n-k))$. By Corollary \ref{cor:wedge(R_k^*)} of
Proposition \ref{prop:wedge(R_k^*)}, it follows that
$(\norm{\cdot}_K^{-k})^\wedge = c(n,k) R_{k}^*(d\mu^\perp)$,
implying that $(\norm{\cdot}_K^{-k})^\wedge \geq 0$, and hence $K
\in \I_k^n$. By Lemma \ref{lem:BP-dense}, and the fact that
$\I_k^n$ is closed under limit in the radial metric, this
concludes the proof.
\end{proof}

Applying Proposition \ref{prop:wedge(R_k^*)} to the function
$f=0$, once for $m=k$ and once for $m=n-k$, we also immediately
deduce the following useful:
\begin{prop} \label{prop:same-kernels}
\[
Ker R_{n-k}^* = Ker R_k^* \circ I.
\]
\end{prop}
This is equivalent by a standard duality argument to the
following Proposition, which may be deduced directly from Theorem
\ref{thm:R_k(wedge)}:
\begin{prop} \label{prop:same-images}
\[
\overline{Im R_{n-k}} = \overline{Im I \circ R_k}.
\]
\end{prop}

\medskip

We conclude this section by introducing a family of very natural
operators acting on $C(G(n,k))$ to itself, and showing a few nice
properties which they share. Denote by $V_k: C(G(n,k)) \rightarrow
C(G(n,k))$ the operator defined as $V_k = I \circ R_{n-k} \circ
R_{k}^*$.

\begin{prop}
$V_k$ is self-adjoint.
\end{prop}
\begin{proof}
It is actually not hard to show this directly, just by using
double-integration as in Section \ref{sec:2}. Nevertheless, we
prefer to use Proposition \ref{prop:wedge(R_k^*)}. Let $f,g \in
C^\infty(G(n,n-k))$. Then by Proposition \ref{prop:wedge(R_k^*)},
the Spherical Parseval identity and Proposition
\ref{prop:wedge(R_k^*)} again, we have:
\begin{eqnarray}
\nonumber \scalar{V_{n-k}(f),g}_{G(n,n-k)} &=&
\scalar{R_{n-k}^*(f),(I\circ R_k)^*(g)} \\
\nonumber &=& c(n,k)^{-1} \scalar{R_{n-k}^*(f),(E_{-k}^\wedge
\circ
R_{n-k}^*)(g)} \\
\nonumber &=& c(n,k)^{-1} \scalar{(E_{-k}^\wedge \circ
R_{n-k}^*)(f),R_{n-k}^*(g)} \\
\nonumber &=& \scalar{(I\circ R_k)^*(f),R_{n-k}^*(g)} =
\scalar{f,V_{n-k}(g)}_{G(n,n-k)}
\end{eqnarray}
Since $C^\infty(G(n,n-k))$ is dense in $C(G(n,n-k))$ in the
maximum norm, and the operators $R_{n-k}^*$ and $R_k$, and hence
$V_{n-k}$, are continuous w.r.t. this norm, it follows that the
same holds for any $f,g \in C(G(n,n-k))$.
\end{proof}

\begin{prop}
\[
V_{n-k} = I \circ V_{k} \circ I.
\]
\end{prop}
\begin{proof}
This time we give the proof in operator style notations. The
formal details are filled in exactly the same manner as above.
Using the definition of $V_k$, and the identities
(\ref{eq:wedge(R_k^*)}) and (\ref{eq:R_k(wedge)-operator}), we
have:
\begin{eqnarray}
\nonumber I \circ V_{k} \circ I & = & R_{n-k} \circ R_k^* \circ
I =  c(n,k)^{-1} R_{n-k} \circ E_{-k}^\wedge \circ R_{n-k}^* \\
\nonumber & = & I \circ R_k \circ R_{n-k}^* = V_{n-k}.
\end{eqnarray}
\end{proof}

It is known (e.g. \cite{Gelfand-Graev-Rosu}) that for $1<k<n-1$,
even if we restrict the operators $R_m$ to infinitely smooth
functions, $Ker R_k^{*} \neq \set{0}$ and $\overline{Im R_{n-k}}
\neq C^\infty(G(n,n-k))$, and therefore $V_k$ is neither
injective nor surjective onto a dense set for those values of
$k$. Since $\overline{Im R_k^{*}} = C_e(S^{n-1})$ and $Ker
R_{n-k} = \set{0}$, it follows that $Ker V_k = Ker R^*_k$ and
$\overline{Im V_k} = \overline{Im I \circ R_{n-k}} = \overline{Im
R_k}$ (by Proposition \ref{prop:same-images}). A standard duality
argument shows that $\overline{Im R_k}$ is orthogonal to $Ker
R^*_k$ (as measures acting on continuous functions, and therefore
as functions when $R^*_k$ is restricted to $C(G(n,k))$), and
therefore we may consider $V_k$ as an operator from $\overline{Im
R_k}$ to $\overline{Im R_k}$, which is injective and surjective
onto a dense set. A natural question for Integral Geometrists
would be to find a nice inversion formula for $V_k$. Note that by
a standard double-integral argument, the operator $R_k^* \circ I
\circ R_{n-k} : S^{n-1} \rightarrow S^{n-1}$ is exactly the usual
Spherical Radon transform $R$ (for every $k$), and under the
standard identification between $G(n,n-1)$, $G(n,1)$ and
$S^{n-1}$, so are $V_1$ and $V_{n-1}$.

%%%%%%%%%%%%%%%%%%%%%%%%%%%%%%%%%%%%%%%%%%%%%%%%%%%%%%%%%%%%%%%%%%%%%%%%%%%%%%%%%%%%%%%
%%%%%%%%%%%%%%%%%%%%%%%%%%%%%%%%%%%%%%%%%%%%%%%%%%%%%%%%%%%%%%%%%%%%%%%%%%%%%%%%%%%%%%%
%%%%%%%%%%%%%%%%%%%%%%%%%%%%%%%%%%%%%%%%%%%%%%%%%%%%%%%%%%%%%%%%%%%%%%%%%%%%%%%%%%%%%%%
%%%%%%%%%%%%%%%%%%%%%%%%%%%%%%%%%%%%%%%%%%%%%%%%%%%%%%%%%%%%%%%%%%%%%%%%%%%%%%%%%%%%%%%
%%%%%%%%%%%%%%%%%%%%%%%%%%%%%%%%%%%%%%%%%%%%%%%%%%%%%%%%%%%%%%%%%%%%%%%%%%%%%%%%%%%%%%%
%%%%%%%%%%%%%%%%%%%%%%%%%%%%%%%%%%%%%%%%%%%%%%%%%%%%%%%%%%%%%%%%%%%%%%%%%%%%%%%%%%%%%%%
%%%%%%%%%%%%%%%%%%%%%%%%%%%%%%%%%%%%%%%%%%%%%%%%%%%%%%%%%%%%%%%%%%%%%%%%%%%%%%%%%%%%%%%
%%%%%%%%%%%%%%%%%%%%%%%%%%%%%%%%%%%%%%%%%%%%%%%%%%%%%%%%%%%%%%%%%%%%%%%%%%%%%%%%%%%%%%%

\section{Equivalent formulations of $\BP_k^n = \I_k^n$}
\label{sec:4}

In this section we use the results and techniques of the previous
sections together with basic tools from Functional Analysis to
derive equivalent formulations of the natural conjecture that
$\BP_k^n = \I_k^n$. As mentioned in the Introduction, the
relevance of this conjecture to Convex Geometry stems from the
generalized $k$-codimensional Busemann-Petty problem. It was
shown in \cite{Zhang-Gen-BP} that the answer to this problem is
positive iff every convex body in $\Real^n$ is in $\BP_k^n$, and
this was shown to be false
(\cite{Bourgain-Zhang},\cite{Koldobsky-I-equal-BP}) for $k<n-3$,
but the cases of $k=n-3$ and $k=n-2$ remain open. The analogous
question for $\I_k^n$ turned out to be easier using the analytic
tools provided by the Fourier transform, and it was shown by
Koldobsky in \cite{Koldobsky-convex-is-n-3-intersection} that
$\I_k^n$ contains all $n$-dimensional convex bodies iff $k\geq
n-3$. Hence, a positive answer to whether $\BP_k^n = \I_k^n$ would
imply a positive answer to the generalized $k$-codimensional
Busemann-Petty problem, for $k \geq n-3$. The equivalent
formulations derived in this section indicate that the $\BP_k^n =
\I_k^n$ question is connected and equivalent to very natural
questions in Integral Geometry.

Before we start, we would like to give an intuitive equivalent
formulation to $\BP_k^n = \I_k^n$. By Grinberg and Zhang's
characterization (Theorem \ref{thm:G&Z}), $\BP_k^n$ is exactly the
class of star-bodies generated from the Euclidean Ball $D_n$ by
means of full-rank linear transformations, $k$-radial sums, and
limit in the radial metric. Loosely speaking, we say that "modulo
these operations", $D_n$ is the only member of $\BP_k^n$. Since
$\I_k^n$ is closed under these operations as well, we can ask
whether "modulo these operations" $D_n$ is the only star-body such
$(\norm{\cdot}_{D_n}^{-k})^\wedge \geq 0$. In terms of functions
on the sphere, this is equivalent to asking whether "modulo these
operations", the only function $f \in C^\infty_e(S^{n-1})$ such
that $f \geq 0$ and $E_{-k}^\wedge(f) \geq 0$ is the constant
function $f=1$ (note that we may restrict our attention to
infinitely smooth functions because of Lemma \ref{lem:I-dense}).
This formulation transforms the problem to the language of Fourier
transforms. As opposed to this, our other formulations in this
section will use the language of the Radon transforms and Integral
Geometry.

\medskip

We will use the following notations. $R_{m}(C(S^{n-1}))_{+}$ will
denote the non-negative functions in the image of $R_{m}$ and
$R_{m}(C_+(S^{n-1}))$ will denote the image of $R_{m}$ acting on
the cone $C_+(S^{n-1})$ (which is the same as its image acting on
$C_{+,e}(S^{n-1})$).% and $\overline{A}$ will denote the closure of
%$A$ w.r.t. the underlying norm.
We denote $G = G(n,n-k)$ for short.

\smallskip

It is well known (e.g.
\cite{Gelfand-Graev-Rosu},\cite{Helgason-Book},\cite{Strichartz})
that $R_{n-k}:C_e(S^{n-1}) \rightarrow C(G(n,n-k))$ is an
injective operator, but it is not onto for $k<n-1$, and
$\overline{Im R_{n-k}} \neq C(G(n,n-k))$ for $1<k<n-1$. We will
restrict our discussion to this range of $k$. It follows by an
elementary duality argument, that the image of the dual operator
$R_{n-k}^*: \M(G(n,n-k)) \rightarrow \M_e(S^{n-1})$ is dense in
$\M_e(S^{n-1})$ in the $w^*$-topology, but $R_{n-k}^*$ is not
injective and has a non-trivial kernel. It is known that the
dense image in $\M_e(S^{n-1})$ contains $C_e^\infty(S^{n-1})$,
and in fact an explicit inversion formula was obtained by
Koldobsky in \cite[Proposition 3]{Koldobsky-I-equal-BP} (which is
not unique because of the kernel). It follows from Koldobsky's
argument (or from the general results of
\cite{Gelfand-Graev-Rosu}) that:
\begin{lem} \label{lem:smooth-pre-image}
If $f \in C_e^\infty(S^{n-1})$ then there exists a $g \in
C^\infty(G(n,n-k))$ such that $f = R_{n-k}^*(g)$.
\end{lem}
\noindent It will also be useful to note that:
\begin{equation} \label{eq:KerR}
Ker R_{n-k}^* = \set{\mu \in \M(G(n,n-k)) | \scalar{\mu,f} = 0
\;\; \forall f \in Im R_{n-k}},
\end{equation}
and to recall Propositions \ref{prop:same-kernels} and
\ref{prop:same-images}, which show that $Ker R_{n-k}^* = Ker
R_k^* \circ I$ and $\overline{Im R_{n-k}} = \overline{Im I \circ
R_k}$. The latter immediately implies:
\begin{equation} \label{eq:positive-images}
 \overline{R_{n-k}(C_{+}(S^{n-1}))} , \overline{I \circ R_k
(C_{+}(S^{n-1}))} \subset \overline{R_{n-k}(C(S^{n-1}))_{+}}.
\end{equation}

\medskip

It will be useful to consider the quotient space:
\[\M(n,n-k) = \M(G(n,n-k)) / Ker R_{n-k}^*,\]
which is the space of bounded linear functionals on the subspace
$\overline{Im R_{n-k}}$ of $C(G(n,n-k))$. By abuse of notation,
we will also think of $R_{n-k}^*$ as an operator from $\M(n,n-k)$
to $\M_e(S^{n-1})$, and although this does not change its image,
it is now injective on $\M(n,n-k)$. The same is true for $R_k^*
\circ I$, since $Ker R_{n-k}^* = Ker R_k^* \circ I$, and we may
proceed to interpret $R_{n-k}^*(d\mu)$ and $R_k^*(d\mu^\perp)$ in
the usual way for $\mu \in\M(n,n-k)$, since these values are the
same for the entire co-set $\mu + Ker R_{n-k}^*$. If $R_{n-k}^*$
were onto $\M_e(S^{n-1})$, or even $C_e(S^{n-1})$, we could
proceed by identifying between a star-body $K$ and a signed Borel
measure $\mu$ in $\M(n,n-k))$, by the correspondence
$\norm{\cdot}_K^{-k} = R_{n-k}^*(d\mu)$. Unfortunately, the
general theory does not guarantee this, and in fact we believe
that some star-bodies do not admit such a representation
(although we have not been able to find a reference for this).
%Unfortunately this is not so, as some star-bodies do not admit
%such a representation (although we have not been able to find a
%reference for this).
But as remarked earlier, $C_e^\infty(S^{n-1})$ does lie in the
image of $R_{n-k}^*$, and this is enough for our purposes.

Let us now review the definitions of $\BP_k^n$ and $\I_k^n$. Our
original definition required that $K \in \BP_k^n$ iff $\rho_K^k =
R_{n-k}^*(d\mu)$ for some non-negative measure $\mu \in
\M_+(G(n,n-k))$. We claim that this is equivalent to requiring
that $\mu \in \M_+(n,n-k)$, since by a version of the Hahn-Banach
Theorem (\cite[Lemma 4.3]{Grinberg-Zhang}), any non-negative
functional on $\overline{Im R_{n-k}}$ may be extended to a
non-negative functional on the entire $C(G(n,n-k))$, and the
converse is trivially true. Defining $\M(\BP_k^n)$ as the set of
non-negative functionals in $\M(n,n-k)$:
\[
\M(\BP_k^n) = \M_+(n,n-k),
\]
we see that:
\begin{lem} \label{lem:alternative-defn-BP}
Let $K$ be a star-body in $\Real^n$. Then $K \in \BP_k^n$ iff
$\rho_K^k = R_{n-k}^*(d\mu)$, for some $\mu \in \M(\BP_k^n)$.
\end{lem}

\noindent Let us also define $\M(\I_k^n)$ as:
\[
\M(\I_k^n) = \set{ \mu \in \M(n,n-k) \; | \; R_{n-k}^*(d\mu) \geq
0 \;,\; R_{k}^*(d\mu^\perp) \geq 0},
\]
where ``$\nu \geq 0$" means that $\nu$ is a non-negative measure
in $\M_e(S^{n-1})$. Using co-set notations, let us also define:
\[
\M^\infty(n,n-k) = \set{f + Ker R_{n-k}^* \; | \; f \in
C^\infty(G) },
\]
and denote:
\[
\M^\infty(\I_k^n) = \M(\I_k^n) \cap \M^\infty(n,n-k),
\] and
\[
\M^\infty_+(n,n-k) = \M^\infty(\BP_k^n) = \M(\BP_k^n) \cap
\M^\infty(n,n-k).
\]

Unfortunately, we cannot give a completely analogous
characterization to Lemma \ref{lem:alternative-defn-BP} for
$\I_k^n$ and $\M(\I_k^n)$. However, we have the following:
\begin{lem} \label{lem:alternative-defn-I}
Let $K$ be an \emph{infinitely smooth} star-body in $\Real^n$.
Then $K \in \I_k^n$ iff $\rho_K^k = R_{n-k}^*(d\mu)$, for some
$\mu \in \M^\infty(\I_k^n)$.
\end{lem}
\begin{proof}
We will first prove the "only if" part. Assume that $K \in
\I_k^{n,\infty}$. By Lemma \ref{lem:smooth-pre-image}, there
exists a signed measure $\mu \in \M^\infty(n,n-k)$ so that
$\norm{\cdot}_K^{-k} = R_{n-k}^*(d\mu)$. By Corollary
\ref{cor:wedge(R_k^*)} of Proposition \ref{prop:wedge(R_k^*)}, it
follows that $(\norm{\cdot}_K^{-k})^\wedge = c(n,k)
R_k^*(d\mu^\perp)$. Since $\norm{\cdot}_K^{-k} \geq 0$ because
$K$ is a star-body and $(\norm{\cdot}_K^{-k})^\wedge \geq 0$
because $K \in \I_k^n$, it follows that $R_{n-k}^*(d\mu)\geq 0$
and $R_k^*(d\mu^\perp) \geq 0$, proving that $\mu \in
\M^\infty(\I_k^n)$. The "if" part follows from Corollary
\ref{cor:wedge(R_k^*)} in exactly the same manner, since
$(\norm{\cdot}_K^{-k})^\wedge = c(n,k) R_k^*(d\mu^\perp) \geq 0$
for a measure $\mu \in \M^\infty(\I_k^n)$ such that
$\norm{\cdot}_K^{-k} = R_{n-k}^*(d\mu)$.
\end{proof}
\begin{rem}
It seems that any attempt to prove the "only if" part of the lemma
for a general star-body $K \in \I_k^n$ by approximating it with
$K_i \in \I_k^{n,\infty}$ will fail. The reason is that we have
no way of controlling the norm of the (a-priori signed) measures
$\mu_i \in \M(\I_k^n)$ for which $\rho_{K_i}^k =
R_{n-k}^*(d\mu_i)$, and therefore it is not guaranteed that
$\mu_i$ will converge to some measure (like in the usual argument
which uses the $w^*$-compactness of the unit-ball of
$\M(n,n-k)$). If it were known that the $\mu_i$ are non-negative
(this would follow if $\M(\BP_k^n) = \M(\I_k^n)$), it would
follow that $\norm{\mu_i} = \norm{R_{n-k}^*(d\mu_i)}$ (since
$R_{n-k}^*(d\mu_i)$ is non-negative), and over the latter term we
do have control. The "if" part of the lemma may be proved without
any smoothness assumption by the standard approximation argument.
\end{rem}

We now see that we have derived alternative definitions of
$\BP_k^n$ and $\I_k^{n,\infty}$ using a common language of Radon
transforms and without using the Fourier transform. Note that
even if we could remove the restriction of infinite smoothness
from Lemma \ref{lem:alternative-defn-I}, it would not be yet clear
that $\BP_k^n = \I_k^n$ iff $\M(\BP_k^n) = \M(\I_k^n)$, since for
a general $\mu \in \M(\BP_k^n)$ or $\mu \in \M(\I_k^n)$,
$R_{n-k}^*(d\mu)$ may not be a measure with continuous density
(and hence cannot equal $\rho_K^k$ for a star-body $K$). We do
however have:
\begin{lem} \label{lem:M(BP)-subset-M(I)}
\[\M(\BP_k^n) \subset \M(\I_k^n) \]
\end{lem}
\begin{proof}
If $\mu \in \M_+(n,n-k)$ then trivially $R_{n-k}^*(d\mu)\geq 0$
and $R_{k}^*(d\mu^\perp)\geq 0$, hence $\mu \in \M(\I_k^n)$.
Although the proof is trivial, note that underlying this
statement are Propositions \ref{prop:same-kernels} and
\ref{prop:same-images} which enabled us to restrict $R_{n-k}^*$
and $R_k^* \circ I$ to $\M(n,n-k)$.
\end{proof}

We may now formulate the main Theorem of this section:

\begin{thm} \label{thm:equivalence}
Let $n$ and $1 \leq k \leq n-1$ be fixed. Then the following are
equivalent:
\begin{enumerate}
\item
\[
\BP_k^n = \I_k^n.
\]
\item
\[
\M^\infty(\BP_k^n) = \M^\infty(\I_k^n).
\]
\item
\[
\M(\BP_k^n) = \M(\I_k^n).
\]
\item
\[
\overline{R_{n-k}(C(S^{n-1}))_{+}} =
\overline{R_{n-k}(C_{+}(S^{n-1})) + I \circ R_k (C_{+}(S^{n-1}))}.
\]
\item
If $\mu + 1 \in \M(\BP_k^n)$ and $\mu \in \M(\I_k^n)$, then $\mu
\in \M(\BP_k^n)$.
\item
There does not exist a measure $\mu \in \M^\infty_+(n,n-k)$ such
that $R_{n-k}^*(d\mu) \geq 1$ and $R_k^*(d\mu^\perp) \geq 1$
(where ``$\nu \geq 1$" means that $\nu-1$ is a non-negative
measure), and such that:
\begin{equation} \label{eq:infimum}
\inf \set{\scalar{\mu,f} | f \in R_{n-k}(C(S^{n-1}))_{+} \textrm{
and } \scalar{1,f} = 1 } = 0.
\end{equation}
\end{enumerate}
\end{thm}

We will show $(2) \Rightarrow (1)$, $(1) \Rightarrow (3)$, $(3)
\Leftrightarrow (4)$, $(5) \Rightarrow (6)$ and $(6) \Rightarrow
(2)$. Obviously, $(3) \Rightarrow (2)$ and $(3) \Rightarrow (5)$.

\begin{proof}[Proof of $(2) \Rightarrow (1)$]
Let $K \in \I_k^{n,\infty}$. In view of Lemma \ref{lem:I-dense},
the fact that $\BP_k^n$ is closed under limit in the radial
metric, and Corollary \ref{cor:BP-subset-I}, it is enough to show
that $K \in \BP_k^n$. By Lemma \ref{lem:alternative-defn-I},
$\rho_K^k = R_{n-k}^*(d\mu)$ for some $\mu \in
\M^\infty(\I_k^n)$. By our assumption that $\M^\infty(\BP_k^n) =
\M^\infty(\I_k^n)$ and by Lemma \ref{lem:alternative-defn-BP}, it
follows that $K \in \BP_k^n$ (in fact $K \in \BP_k^{n,\infty}$).
\end{proof}

\begin{proof}[Proof of $(1) \Rightarrow (3)$]
In view of Lemma \ref{lem:M(BP)-subset-M(I)}, it is enough to
prove $\M(\I_k^n) \subset \M(\BP_k^n)$. Let $\mu \in \M(\I_k^n)$,
so $R_{n-k}^*(d\mu) \geq 0$ and $R_k^*(d\mu^\perp) \geq 0$. Let
$\set{u_i} \subset C^\infty(O(n))$ be an approximate identity as
in Lemma \ref{lem:approximate-identity}. Let $K_i$ denote the
infinitely smooth star-body defined by:
\[
\norm{\cdot}_{K_i}^{-k} = u_i \ast R_{n-k}^*(\mu) \geq 0
\]
(we used $R_{n-k}^*(\mu) \geq 0$ to verify that $K_i$ is indeed a
star-body). As in the proof of Lemma \ref{lem:BP-dense}, it is
easy to see that $\norm{\cdot}_{K_i}^{-k} = R_{n-k}^*(u_i \ast
\mu)$, so by Corollary \ref{cor:wedge(R_k^*)} of Proposition
\ref{prop:wedge(R_k^*)} we have:
\[
(\norm{\cdot}_{K_i}^{-k})^\wedge = c(n,k) R_{k}^*((u_i \ast
\mu)^\perp) = R_{k}^*(u_i \ast \mu^\perp) = u_i \ast
R_{k}^*(\mu^\perp) \geq 0.
\]
Hence $K_i \in \I_k^n$, and by our assumption that $\BP_k^n =
\I_k^n$, it follows that $K_i \in \BP_k^n$. By Lemma
\ref{lem:alternative-defn-BP}, this implies that
$\norm{\cdot}_{K_i}^{-k} = R_{n-k}^*(d\eta_i)$, where $\eta_i \in
\M(BP_k^n)$. The injectivity of $R_{n-k}^*$ on $\M(n,n-k)$
implies that $u_i \ast \mu = \eta_i \in \M(BP_k^n)$. Lemma
\ref{lem:approximate-identity} shows that $u_i \ast \mu$ tends to
$\mu$ in the $w^*$-topology, and since $\M(BP_k^n)$ is obviously
closed in this topology, it follows that $\mu \in \M(BP_k^n)$.
\end{proof}

For the proof of $(3) \Leftrightarrow (4)$ and for later use, we
will need to recall a few classical notions from Functional
Analysis (e.g. \cite{Bourbaki-TVS}). A cone $P$ in a Banach space
$X$ is a non-empty subset of $X$ such that $x,y \in P$ implies
$c_1 x + c_2 y \in P$ for every $c_1,c_2 \geq 0$. The dual cone
$P^* \subset X^*$ is defined by $P^* = \set{ x^* \in X^* |
\scalar{x^*,p} \geq 0 \;\; \forall p\in P}$. Therefore $P^*$ is
always closed in the $w^*$-topology, and $P^* =
(\overline{P})^*$. It is also easy to check that $P_1 \subset
P_2$ implies $P_2^* \subset P_1^*$, $(P_1 + P_2)^* = P_1^* \cap
P_2^*$ and $(P_1 \cap P_2)^* = P_1^* + P_2^*$. An immediate
consequence of the Hahn-Banach Theorem is that $\overline{P_1} =
\overline{P_2}$ iff $P_1^* = P_2^*$.

\begin{proof}[Proof of  $(3) \Leftrightarrow (4)$]
All the sets appearing in (3) and (4) are clearly cones. It
remains to show that the cones in both sides of (3) are exactly
the dual cones to the ones in both sides of (4). The equivalence
then follows by the Hahn-Banach Theorem, as in the last statement
of the previous paragraph.

By definition, $\M(\BP_k^n)$ is dual to
$\overline{R_{n-k}(C(S^{n-1}))_{+}}$. The cones:
\[\set{ \mu \in
\M(n,n-k) \; | \; R_{n-k}^*(d\mu) \geq 0 } \] and \[\set{ \mu \in
\M(n,n-k) \; | \; R_{k}^*(d\mu^\perp) \geq 0 }\] are immediately
seen to be dual to $R_{n-k}(C_{+}(S^{n-1}))$ and $I \circ R_k
(C_{+}(S^{n-1}))$, respectively. Since $(P_1 + P_2)^* = P_1^*
\cap P_2^*$, it follows that:
\[
\M(\I_k^n) = \brac{\overline{R_{n-k}(C_{+}(S^{n-1})) + I \circ
R_k (C_{+}(S^{n-1}))}}^*.
\]
This concludes the proof.
\end{proof}
\begin{rem}
By (\ref{eq:positive-images}), we have:
\[
\overline{R_{n-k}(C(S^{n-1}))_{+}} \supset
\overline{R_{n-k}(C_{+}(S^{n-1})) + I \circ R_k (C_{+}(S^{n-1}))}.
\]
By duality, we see again that:
\[
\M(\BP_k^n) \subset \M(\I_k^n).
\]
\end{rem}

\begin{proof}[Proof of $(5) \Rightarrow (6)$]
This follows immediately from the definitions. Assume that (6) is
false, so that there exists a measure $\mu \in
\M_+^\infty(n,n-k)$ such that $R_{n-k}^*(d\mu) \geq 1$ and
$R_{k}^*(d\mu^\perp) \geq 1$ and such that (\ref{eq:infimum})
holds. Define $\mu' = \mu - 1$, and so $\mu'+1 \in \M(\BP_k^n)$,
$\mu' \in \M(\I_k^n)$, and (\ref{eq:infimum}) shows that $\mu'$
is not in $\M(\BP_k^n)$. Therefore $\mu'$ is a counterexample to
(5).
\end{proof}

\begin{proof}[Proof of $(6) \Rightarrow (2)$]
Assume that $(2)$ is false, so $\M^\infty(\BP_k^n) \neq
\M^\infty(\I_k^n)$. By Lemma \ref{lem:M(BP)-subset-M(I)}, this
means that there exists a measure $\mu' \in \M^\infty(\I_k^n)
\setminus \M(\BP_k^n)$. Since $\mu' \in \M^\infty(n,n-k)$, we can
write $\mu' = g + Ker R_{n-k}^*$ with $g \in C^\infty(G(n,n-k))$.
Assume that $\min(g) = -C$ where $C>0$, otherwise we would have
$\mu' \in \M(\BP_k^n)$.

Now consider the measure $\mu_\lambda = (1-\lambda) \mu' +
\lambda \in \M^\infty(n,n-k)$ for $\lambda \in [0,1]$. Since
$\M(\BP_k^n)$ is convex, contains the measure $1$, and is closed
in the $w^*$-topology, it follows that there exists a $\lambda_0
\in (0,1]$ so that $\mu_\lambda \in \M^\infty(\BP_k^n)$ iff
$\lambda \in [\lambda_0,1]$. But for $\lambda_1 = C / (1+C)$ we
already see that $\mu_{\lambda_1} \in \M(\BP_k^n)$, because
$\mu_{\lambda_1} = g_{\lambda_1} + Ker R_{n-k}^*$ and
$g_{\lambda_1} = 1/(1+C) g + 1-1/(1+C) \in C^\infty_+(G(n,n-k))$.
We conclude that $\lambda_0 \in (0,1)$.

Now define $\mu = \mu_{\lambda_0} / \lambda_0 \in \M(\BP_k^n)$,
and notice that $\mu - 1 = (1-\lambda_0)/\lambda_0 \mu' \in
\M^\infty(\I_k^n)$, implying that $R_{n-k}^*(d\mu) \geq 1$ and
$R_k^*(d\mu^\perp) \geq 1$. It remains to show (\ref{eq:infimum}).
Assume by negation that:
\[
\inf \set{\scalar{\mu,f} | f \in R_{n-k}(C(S^{n-1}))_{+} \textrm{
and } \scalar{1,f} = 1 } = \delta > 0.
\]
But then it is easy to check that for $\lambda_2 = \lambda_0
(1-\delta)/ (1 - \delta \lambda_0) < \lambda_0$,
$\scalar{\mu_{\lambda_2},f} \geq 0$ for all $f \in
R_{n-k}(C(S^{n-1}))_{+}$, and hence for all $f \in
\overline{R_{n-k}(C(S^{n-1}))_{+}}$. Therefore $\mu_{\lambda_2}
\in \M^\infty(\BP_k^n)$, in contradiction to the definition of
$\lambda_0$. Therefore (\ref{eq:infimum}) is shown, concluding
the proof.

\end{proof}

\begin{rem} \label{rem:back-to-M}
In formulation $(6)$, it is equivalent to require that $\mu \in
\M_+(n,n-k)$ and also $\mu \in \M(G(n,n-k))$ instead of $\mu \in
\M_+^\infty(n,n-k)$. The equivalence of $\mu \in \M_+(n,n-k)$
follows since we have not used the fact that $\mu \in
\M^\infty(n,n-k)$ in the proof (by negation) of $(5) \Rightarrow
(6)$. The equivalence of $\mu \in \M(G(n,n-k))$ follows by the
previously mentioned version of the Hahn-Banach Theorem (which
was used to derive Lemma \ref{lem:alternative-defn-BP}). This is
the formulation which was used in the Introduction.
\end{rem}

\medskip

We proceed to develop several more formulations of the $\BP_k^n =
\I_k^n$ question. Unfortunately, we cannot show an equivalence
with the original question, but rather a weak type of implication.
We formulate a very natural conjecture, and show that together
with a positive answer to the $\BP_k^n = \I_k^n$ question, the new
formulations are implied.

%Unfortunately, we can only show that these new formulations
%\emph{imply} a positive answer to the latter, and are unable to
%complete the equivalence. However, we formulate a very natural
%conjecture, under which the equivalence is established.

Given an Borel set $Z \subset G(n,n-k)$, we define the restriction
of a measure $\mu \in \M(G(n,n-k))$ to $Z$, denoted $\mu|_Z \in
\M(G(n,n-k))$, as the measure satisfying $\mu|_Z(A) = \mu(A \cap
Z)$ for any Borel set $A \subset G(n,n-k)$. We will say that
$\mu$ is supported in a closed set $Z$, if $\mu|_{Z^C} = 0$, and
define the support of $\mu$, denoted $supp(\mu)$, as the minimal
closed set $Z$ in which $\mu$ is supported (it is easy to check
that this is well-defined). It is also easy to check that:
\begin{lem} \label{lem:supp-and-zero-set}
If $f \in C(G(n,n-k))$, $\mu \in \M(G(n,n-k))$ and $supp(\mu)
\subset f^{-1}(0)$ then $\scalar{\mu,f} = 0$. Conversely, if $f
\in C_+(G(n,n-k))$, $\mu \in \M_+(G(n,n-k))$ and $\scalar{\mu,f} =
0$, then $supp(\mu) \subset f^{-1}(0)$.
\end{lem}
We also recall the definition of the Covering Property from the
Introduction. A set closed set $Z \subset G(n,n-k)$ is said to
satisfy the \emph{covering property} if:
%\smallskip
%\noindent\textbf{The Covering Property. }
\begin{equation} \label{eq:condition-on-A-again}
\bigcup_{E \in Z} E \cap S^{n-1} = S^{n-1} \; \text{ and } \;
\bigcup_{E \in Z} E^\perp \cap S^{n-1} = S^{n-1}.
\end{equation}

Our starting point is formulation $(6)$ in Theorem
\ref{thm:equivalence}, which involves both a function $f$ and a
measure $\mu$. Note that the requirement that if $f \in
\overline{R_{n-k}(C(S^{n-1}))_{+}}$ and $\scalar{1,f} = 1$, then
$\scalar{\mu,f}$ is bounded away from zero, is stronger than
demanding that $\scalar{\mu,f} \neq 0$. The motivation for the
following discussion stems from the impression that the
conditions on $\mu$, namely that $\mu \in \M_+(G(n,n-k))$
(following Remark \ref{rem:back-to-M}), $R_{n-k}^*(d\mu) \geq 1$
and $R_{k}^*(d\mu^\perp) \geq 1$, may be equivalently specified
by some condition on the support of $\mu$. In that case, the
condition that $\scalar{\mu,f} \neq 0$ becomes a condition on the
set $f^{-1}(0)$. Let us show the following necessary condition on
the support of such a $\mu$ as above:

\begin{lem} \label{lem:covering-necessary}
Let $\mu \in \M_+(G(n,n-k))$ so that $R_{n-k}^*(d\mu) \geq 1$ and
$R_{k}^*(d\mu^\perp) \geq 1$. Then $supp(\mu)$ satisfies the
covering property.
\end{lem}

\begin{proof}
Denote by $Z =supp(\mu)$ and $\widetilde{Z} = \bigcup_{E \in Z} E
\cap S^{n-1}$. We will show that if $\mu \in \M_+(G(n,n-k))$ and
$R_{n-k}^*(d\mu) \geq 1$ then  $\widetilde{Z} = S^{n-1}$. The
other "half" of the covering property follows similarly from
$R_{k}^*(d\mu^\perp) \geq 1$.

Notice that for $E_1,E_2 \in G(n,n-k)$, the Hausdorff distance
between $E_1 \cap S^{n-1}$ and $E_2 \cap S^{n-1}$ is equivalent
to the distance between $E_1$ and $E_2$ in $G(n,n-k)$. It follows
that since $Z$ is closed, so is $\widetilde{Z}$. Now assume that
$\widetilde{Z} \neq S^{n-1}$, so there exists a $\theta \in
S^{n-1}$ and an $\epsilon>0$, so that $\widetilde{B} =
B_{S^{n-1}}(\theta,\epsilon) \cup B_{S^{n-1}}(-\theta,\epsilon)
\subset \widetilde{Z}^C$. Let $f \in C_{e,+}(S^{n-1})$ be any
non-zero function supported in $\widetilde{B}$. Since
$\widetilde{B} \subset \widetilde{Z}^C$ it follows that $B =
supp(R_{n-k}(f)) \subset Z^C$, and therefore:
\[
\scalar{R_{n-k}^*(\mu),f} = \scalar{\mu,R_{n-k}(f)} = 0.
\]
But on the other hand, since $R_{n-k}^*(d\mu) \geq 1$ and $f\in
C_{e,+}(S^{n-1})$ is non-zero:
\[
\scalar{R_{n-k}^*(\mu),f} \geq \scalar{1,f} > 0,
\]
a contradiction.
\end{proof}

We conjecture that the covering property is also a sufficient
condition in the following sense:

\smallskip
\noindent\textbf{Covering Property Conjecture. } \emph{For any
$n>0$, $1\leq k \leq n-1$, if $Z \subset G(n,n-k)$ is a closed set
satisfying $\bigcup_{E \in Z} E \cap S^{n-1} = S^{n-1}$, then
there exists a measure $\mu \in \M_+(G(n,n-k))$ supported in $Z$,
such that $R_{n-k}^*(d\mu) \geq 1$.}
\smallskip

Under this conjecture, we immediately have the following
counterpart to Lemma \ref{lem:covering-necessary}:

\begin{lem} \label{lem:covering-sufficient}
Assume the Covering Property Conjecture, and let $Z \subset
G(n,n-k)$ be a closed set satisfying the covering property. Then
there exists a measure $\mu \in \M_+(G(n,n-k))$ supported in $Z$,
such that $R_{n-k}^*(d\mu) \geq 1$ and $R_{k}^*(d\mu^\perp) \geq
1$.
\end{lem}
\begin{proof}
Apply the Conjecture to the closed sets $Z \subset G(n,n-k)$ and
$Z^\perp \subset G(n,k)$, and let $\mu_1 \in \M_+(G(n,n-k))$ and
$\mu_2 \in \M_+(G(n,k))$ be the resulting measures. Then $\mu_1 +
\mu_2^\perp$ is supported in $Z$ and satisfies the requirements.
\end{proof}

\begin{rem} \label{rem:prove-covering-conj}
A very natural way to approach the proof of the Covering Property
Conjecture, is to assume that the closed set $Z$ satisfying
$\bigcup_{E \in Z} E \cap S^{n-1} = S^{n-1}$ is minimal w.r.t.
set inclusion (indeed, by Zorn's lemma it is easy to verify that
there exists such a minimal set). The natural candidate for a
measure supported on $Z$ is simply the Hausdorff measure $H_Z$ on
$Z$, and it remains to show that $H_Z$ is a finite measure and
that $R_{n-k}^*(dH_Z) \geq \epsilon$ for some $\epsilon>0$, using
the minimality of $Z$. In particular, one has to show that the
Hausdorff dimension of $Z$ is $k$. Although having some progress
in this direction, we have not been able to give a complete
proof. We also remark that it is easy to construct a
\emph{non-bounded} measure $\mu$ supported on $Z$ for which
$R_{n-k}^*(d\mu) \geq 1$, simply by using the counting measure on
$Z$, i.e. $\mu(A) = \abs{\set{A \cap Z}}$ for any Borel set $A
\subset G(n,n-k)$ (where $\abs{A}$ denotes the cardinality of
$A$).
\end{rem}

As opposed to Theorem \ref{thm:equivalence}, where $R_{n-k}^*$
was treated as an operator on $\M(n,n-k)$, we now go back to the
original definition of $R_{n-k}^*$ as an operator acting on the
entire $\M(G(n,n-k))$. We summarize this in the following lemma,
abbreviating as usual $G = G(n,n-k)$:

\begin{lem} \label{lem:pass-to-whole-space} \hfill
\begin{enumerate}
\item
\[
\M(n,n-k) = \M(G) / Ker R_{n-k}^*.
\]
\item
\[
\M_+(n,n-k) = \set{ \mu + Ker R_{n-k}^* \; | \; \mu \in \M_+(G)}.
\]
\item
\[
\set{\mu \in \M(G) | \scalar{\mu,f} \geq 0 \; \forall f \in
\overline{R_{n-k}(C(S^{n-1}))_{+}}} = \M_+(G) + Ker R_{n-k}^*.
\]
\end{enumerate}
\end{lem}
\begin{proof}
$(1)$ is simply the definition of $\M(n,n-k)$. $(2)$ follows from
$(3)$, since $\M_+(n,n-k)$ is defined as the cone of non-negative
linear functionals on $\overline{Im R_{n-k}}$, and any linear
functional on the subspace may be extended to the entire space,
hence to $\mu \in \M(G)$. $(3)$ was already implicitly used in
the proof of Lemma \ref{lem:alternative-defn-BP}, but we repeat
the argument once more. The right-hand set is clearly a subset of
the left-hand set, since $Ker R_{n-k}^*$ is perpendicular to
$\overline{Im R_{n-k}}$ by (\ref{eq:KerR}). Conversely, any $\mu$
in the left-hand set is a non-negative linear functional on
$\overline{Im R_{n-k}}$, and by a version of the Hahn-Banach
Theorem (as in the proof of Lemma \ref{lem:alternative-defn-BP}),
may be extended to a $\mu' \in \M_+(G)$. Again by (\ref{eq:KerR}),
the difference $\mu' - \mu$ must lie in $Ker R_{n-k}^*$,
concluding the proof.
\end{proof}

We now state several more formulations, which are shown to be
equivalent each to the other. We then show that under the
Covering Property Conjecture, a positive answer to the $\BP_k^n =
\I_k^n$ question would imply these new statements.
%We then show that these formulations for a set $Z$ satisfying the
%covering property imply that $\BP_k^n = \I_k^n$, and that under
%the Covering Property Conjecture, that they are actually
%equivalent to the latter.
For a closed set $Z \subset G(n,n-k)$,
we denote by $\M(Z)$ the set of all measures in $\M(G(n,n-k))$
supported in $Z$.

\begin{thm} \label{thm:equivalence2}
Let $n$ and $1 \leq k \leq n-1$ be fixed, and let $Z \subset
G(n,n-k)$ denote a closed subset.
%which satisfies the covering property.
Then the following are equivalent:
\begin{enumerate}
\item
There does not
%\marginpar{This implies that
%$\overline{R_{n-k}(C(S^{n-1}))_{+}}$ has no extremal rays!!!}
exist a non-zero $f \in \overline{R_{n-k}(C(S^{n-1}))_{+}}$ such
that $Z \subset f^{-1}(0)$.
\item
\[
\overline{R_{n-k}(C(S^{n-1}))_{+}} \cap \set{f \in C(G) \; | \;\;
f|_Z = 0} = \set{0}.
\]
\item
\[
\M_+(G) + Ker R_{n-k}^* + \M(Z) = \M(G).
\]
\item
There exists a measure $\mu \in \M(G)$ such that $R_{n-k}^*(d\mu)
= 0$ and $\mu = \mu_1 + \mu_2$ where $\mu_i \in \M(G)$, $\mu_1
\geq 1$ and $\mu_2$ is supported in $Z$.
\end{enumerate}
\end{thm}

It is clear that $(2)$ is just a convenient reformulation of
$(1)$. We will show that $(2) \Leftrightarrow (3)$ and $(3)
\Leftrightarrow (4)$.

\begin{proof}[Proof of $(2) \Leftrightarrow (3)$]
Again, we use the Hahn-Banach theorem which shows that for cones,
$\overline{P_1} = \overline{P_2}$ iff $P_1^* = P_2^*$. The dual
cone (in $\M(G)$) to $\overline{R_{n-k}(C(S^{n-1}))_{+}}$ is by
definition:
\[
\set{\mu \in \M(G) | \scalar{\mu,f} \geq 0 \; \forall f \in
\overline{R_{n-k}(C(S^{n-1}))_{+}}},
\]
which by Lemma \ref{lem:pass-to-whole-space} is equal to $\M_+(G)
+ Ker R_{n-k}^*$. The dual cone to $C_Z(G) = \set{f \in C(G) \; |
\;\; f|_Z = 0}$ is obviously $\M(Z)$. Indeed, by definition, if
$\mu \in \M(G)$ is not supported in $Z$, there exists a $f \in
C_Z(G)$ such that $\scalar{\mu,f} \neq 0$ (since $Z$ is closed).
Since also $-f \in C_Z(G)$, either $\scalar{\mu,f}$ or
$\scalar{\mu,-f}$ is negative, and therefore $\mu$ cannot be in
the dual cone to $C_Z(G)$. The dual cone to $\set{0}$ is of
course $\M(G)$. Using $(P_1 \cap P_2)^* = P_1^* + P_2^*$, this
concludes the proof.
\end{proof}

\begin{proof}[Proof of $(3) \Rightarrow (4)$]
Apply $(3)$ with the measure $-1 \in \M(G)$ on the right hand
side. Then there exist measures $\nu_1 \in \M_+(G)$, $\nu_2 \in
Ker R_{n-k}^*$ and $\nu_3 \in \M(Z)$, such that $\nu_1 + \nu_2 +
\nu_3 = -1$. Denoting $\mu = -\nu_2$, $\mu_1 = \nu_1 + 1$ and
$\mu_2 = \nu_3$, $(4)$ follows immediately.
\end{proof}

\begin{proof}[Proof of $(4) \Rightarrow (3)$]
$\C(G)$ is dense in $\M(G)$ in the $w^*$-topology, so it is
enough to show that $(4)$ implies $C(G) \subset \M_+(G) + Ker
R_{n-k}^* + \M(Z)$, as the cones on the right hand side are
closed in this topology. Let $g \in C(G)$, so there exists a
constant $C\geq 0$ such that $g + C \geq 0$, and hence $g + C +
Ker R_{n-k}^* \in \M_+(n,n-k)$. By Lemma
\ref{lem:pass-to-whole-space}, this means that $g + C \in \M_+(G)
+ Ker R_{n-k}^*$, and we see that it is enough to show that the
measure $-C$ is in $\M_+(G) + Ker R_{n-k}^* + \M(Z)$. Since all
of the involved sets are cones, it is enough to show the claim
for the measure $-1$. But this follows from formulation $(4)$ in
the same manner is in the previous proof. Indeed, let $\mu =
\mu_1 + \mu_2$ as assured by $(4)$, where $\mu \in Ker
R_{n-k}^*$, $\mu_1-1 \in \M_+(G)$ and $\mu_2 \in \M(Z)$. Then $-1
= (\mu_1 - 1) -\mu + \mu_2 \in \M_+(G) + Ker R_{n-k}^* + \M(Z)$.
This concludes the proof.
\end{proof}

Comparing formulations (6) in Theorem \ref{thm:equivalence} and
(1) in Theorem \ref{thm:equivalence2} for a set $Z$ satisfying
the covering property, and using Lemmas
\ref{lem:covering-necessary} and \ref{lem:covering-sufficient},
the following should now be clear:

\begin{prop} Let $n$ and $1 \leq k \leq n-1$ be fixed.
Assuming the Covering Property Conjecture, if any of the
formulations in Theorem \ref{thm:equivalence} hold, then so do
any of the formulations in Theorem \ref{thm:equivalence2} for
\emph{any} closed $Z \subset G(n,n-k)$ \emph{satisfying the
covering property}.
\end{prop}
\begin{proof}
The statement follows immediately from the remark before the
Proposition, taking into account Remark \ref{rem:back-to-M} and
Lemma \ref{lem:supp-and-zero-set}.
\end{proof}

% I think that a similar claim is not true:
% for an even f:
% R_{n-k}(f) \geq 0 and R_k(f) \geq 0 implies f \geq 0.
% a counterexample should be f = C + E{-k}^\wedge(g), with R_k(g)>= -C * const,
% where the latter function is not positive, achieving a minimum below -C.
% R_{n-k}(f) \geq 0 always, and we should hope that min(R_k(E{-k}^\wedge(g))) > -C.
% This should be possible...

\section{Appendix}

In the Appendix, we formulate and prove Proposition
\ref{prop:Grassmann}, which is an extended version of the
statement from the Introduction and of Corollary
\ref{cor:Grassmann}. We have left the proof of Proposition
\ref{prop:Grassmann} for the Appendix, since the technique
involved differs from those used in the rest of this note.
Although the proposition is of elementary nature and fairly simple
to prove, we have not been able to find a reference to it in the
literature, so we give a self contained proof here. A similar
formulation of the case $k_1,\ldots,k_r=1$ was given by Blashcke
and Petkantschin (see \cite{Santalo-Book},\cite{Miles} for an
easy derivation), and used by Grinberg and Zhang in
\cite{Grinberg-Zhang} to deduce that $\BP_1^n \subset \BP_l^n$
for all $1\leq l \leq n-1$.

We assume some elementary knowledge of exterior products of
differential forms on homogeneous spaces. A rigorous derivation
may be found in \cite{Santalo-Book}, but we recommend the
intuitive exposition in \cite[Sections 2,3]{Miles}. We will also
use the notations from Section \ref{sec:2}.

We will use the following terminology. For a set of $m$ vectors
$\bar{v} = \set{v_1,\ldots,v_m}$ in a Euclidean space $V$, denote
by $Vol_m(\bar{v}) =
det(\set{\scalar{v_i,v_j}}_{i,j=1}^{m})^{1/2}$, which is exactly
the \emph{$m$-dimensional volume of the parallelepiped spanned by
$\bar{v}$}. If $m = \sum_{i=1}^r k_i$, let $U_i$ be a $k_i$
dimensional subspace of $V$. Choose an arbitrary basis $\bar{u^i}
=  \set{u^i_1,\ldots,u^i_{k_i}}$ of $U_i$ such that
$Vol_{k_i}(\bar{u^i}) = 1$, and let $\bar{u} = \cup_{i=1}^r
\bar{u^i}$. Then the \emph{$m$-dimensional volume of the
parallelepiped spanned by unit volume elements of
$U_1,\ldots,U_r$} is defined as $Vol_m(\bar{u})$. It is easy to
verify that this definition indeed does not depend on the basis
$\bar{u^i}$ chosen for $U_i$, as long as $Vol_{k_i}(\bar{u^i}) =
1$ (this will also be clear from the proof of Proposition
\ref{prop:Grassmann}).

\begin{prop} \label{prop:Grassmann}
Let $n>1$ be fixed, let $d$ be an integer between 0 and $n-1$,
and let $D \in G(n,d)$. For $i=1,\ldots,r$, let $k_i \geq 1$
denote integers whose sum $l$ satisfies $l \leq n-d$. For
$a=1,\ldots,n-d$ denote by $G^{a} = G(n,n-a)$, and by $\mu^{a}_D$
the Haar probability measure on $G^a_D$. For $F \in G^l$ and $a =
1,\ldots,l-1$, denote by $\mu^{a}_F$ the Haar probability measure
on $G^{a}_F$. Denote by $\bar{E} = (E_1,\ldots,E_r)$ an ordered
set with $E_i \in G^{k_i}$. Then for any continuous function
$f(\bar{E}) = f(E_1,\ldots,E_r)$ on $G^{k_1} \times \ldots \times
G^{k_r}$:
\begin{eqnarray}
\nonumber & & \!\!\!\!\!\!\!\!\! \int_{E_1 \in G^{k_1}_D} \cdots
\int_{E_r\in G^{k_r}_D}
f(\bar{E}) d\mu^{k_1}_D(E_1)\cdots d\mu^{k_r}_D(E_r) = \\
\nonumber & & \!\!\!\!\!\!\!\!\! \int_{F \in G^l_D} \int_{E_1 \in
G^{k_1}_F} \cdots \int_{E_r \in G^{k_r}_F} f(\bar{E})
\Delta(\bar{E}) d\mu^{k_1}_F(E_1) \cdots d\mu^{k_r}_F(E_r)
d\mu^{l}_D(F),
\end{eqnarray}
where $\Delta(\bar{E}) = C_{n,\set{k_i},l,d}
\Omega(\bar{E})^{n-d-l}$, $C_{n,\set{k_i},l,d}$ is a constant
depending only on $n,\set{k_i},l,d$, and $\Omega(\bar{E})$
denotes the volume of the $l$-dimensional parallelepiped spanned
by unit volume elements of $E_1^\perp,\ldots,E_r^\perp$.
\end{prop}

\begin{rem}
One way to compute the constant $C_{n,\set{k_i},l,d}$ is to use
the function $f = 1$ in Proposition \ref{prop:Grassmann}. Perhaps
a better way is to follow the proof, which gives:
\[
C_{n,\set{k_i},l,d} = \frac{\abs{G(n-d,n-d-l)} \Pi_{i=1}^r
\abs{G(l,l-k_i)}}{ \Pi_{i=1}^r \abs{G(n-d,n-d-k_i)}},
\]
where $\abs{G(a,b)}$ denotes the volume of the Grassmann Manifold
$G(a,b)$, and is given by (\cite{Miles}):
\begin{equation} \label{eq:volume-of-G}
\abs{G(a,b)} = \frac{\abs{S^{a-1}} \cdots \abs{S^{a-b}}}{
\abs{S^{b-1}} \cdots \abs{S^{0}}},
\end{equation}
where $\abs{S^m}$ denotes the volume of the Euclidean unit sphere
$S^m$ of dimension $m$ (and $\abs{S^{0}} = 2$).
\end{rem}

\begin{proof}[Proof of Proposition \ref{prop:Grassmann}]
We will show that the measures $d\mu^{k_1}_D(E_1)\cdots
d\mu^{k_r}_D(E_r)$ and $\Delta(\bar{E}) d\mu^{k_1}_F(E_1) \cdots
d\mu^{k_r}_F(E_r) d\mu^{l}_D(F)$ with $F = \cap_{a=1}^r E_a$
coincide on a set of measure 1 w.r.t. both measures. It is easy
to verify that the set consisting of all $(E_1,\ldots,E_r)$ such
that $dim(\cap_{a=1}^r E_a) = n-l$ satisfies this requirement, and
therefore $F$ above is in $G(n,n-l)$, hence the second measure is
well defined. Indeed, this set is exactly complementary to the
set of all $(E_1,\ldots,E_r)$ such that $\Omega(\bar{E}) = 0$,
which defines a lower dimensional analytic submanifold of
$G^{k_1} \times \ldots \times G^{k_r}$, hence having measure 0
w.r.t. the first (Haar) measure.

If $J \in G(a,c)$, it is well-known (\cite{Miles}) that the volume
element of $G_J(a,b)$ for $b>c$ at $H \in G_J(a,b)$ is given by:
\begin{equation} \label{eq:basic-wedge}
dG_J(a,b)(H) = \bigwedge_{i=c+1}^b \bigwedge_{j=b+1}^{a} w_{i,j}
\;\;,
\end{equation}
where $w_{i,j} = \scalar{e_i,de_j}$, and $\set{e_1,\ldots,e_a}$ is
any orthonormal basis of $\Real^a$ such that $J =
span\set{e_1,\ldots,e_c}$ and $H = span\set{e_1,\ldots,e_b}$.
Indeed, it is easy to verify that this formula does not depend on
the given orthonormal basis satisfying these conditions, by
changing basis and applying a change of variables formula. With
this normalization, the total volume of $G_J(a,b)$ is
$\abs{G(a-c,b-c)}$, as defined in (\ref{eq:volume-of-G})
(\cite{Miles}). Since $d_1 \wedge d_2 = - d_2 \wedge d_1$, the
volume element is signed, corresponding to the assumed
orientation of the element. However, we will henceforth ignore
the orientation and implicitly take the absolute value in all
exterior products, except where it is mentioned otherwise. Note
also that the skew-symmetry implies $d \wedge d = 0$.

Let $\set{f_1,\ldots,f_d}$ be an orthonormal basis of $D$, and
let $\set{f_1,\ldots,f_{n-l}}$ be a completion to an orthonormal
basis of $F$. For $a = 1,\ldots,r$ let
$\set{e^a_{n-l+1},\ldots,e^a_{n-k_a}}$ be an orthonormal basis of
$F^\perp \cap E_a$, and let $\set{e^a_{n-k_a+1},\ldots,e^a_n}$ be
an orthonormal basis of $E_a^\perp$. For every $a$ we define
$e^a_i = f_i$ for $i = 1 , \ldots, n-l$. Then:
\[
d\mu^{k_1}_F(E_1)\cdots d\mu^{k_r}_F(E_r) = C^1_{n,\set{k_i},l,d}
\bigwedge_{a=1}^r \bigwedge_{i=d+1}^{n-k_a}
\bigwedge_{j=n-k_a+1}^{n} w^a_{i,j}\;\;,
\]
where $w^a_{i,j} = \scalar{e^a_i,de^a_j}$ and
$C^1_{n,\set{k_i},l,d} = (\Pi_{i=1}^r \abs{G(n-d,n-d-k_i)})^{-1}$
accounts for the fact that the measure on the left is normalized
to have total mass 1. Notice that by (\ref{eq:basic-wedge}):
\[
\bigwedge_{a=1}^r \bigwedge_{i=n-l+1}^{n-k_a}
\bigwedge_{j=n-k_a+1}^{n} w^a_{i,j} = C^2_{\set{k_i},l}
d\mu^{k_1}_F(E_1) \cdots d\mu^{k_r}_F(E_r),
\]
where $C^2_{\set{k_i},l} = \Pi_{i=1}^r \abs{G(l,l-k_i)}$. It
remains to show that:
\begin{equation} \label{eq:remains}
\bigwedge_{a=1}^r \bigwedge_{i=d+1}^{n-l}
\bigwedge_{j=n-k_a+1}^{n} w^a_{i,j} = C^3_{n,\set{k_i},l,d}
\Delta(\bar{E}) d\mu^l_D(F).
\end{equation}

Now let $\set{g_{n-l+1},\ldots,g_n}$ denote an orthonormal basis
of $F^\perp$, and denote $\lambda^a_{j,v} = \scalar{e^a_j,g_v}$
for $j,v=n-l+1, \ldots, n$. Hence $e^a_j = \sum_{v = n-l+1}^n
\lambda^a_{j,v} g_v$ and $de^a_j = \sum_{v = n-l+1}^n
(d\lambda^a_{j,v} g_v + \lambda^a_{j,v} dg_v)$. Denoting $w_{j,v}
= \scalar{f_j,dg_v}$, we see that since $\scalar{f_i,g_v}=0$, then
for $i=1,\ldots,n-l$ and $j = n-l+1,\ldots,n$:
\begin{equation} \label{eq:w}
w^a_{i,j} = \sum_{v = n-l+1}^n \lambda^a_{j,v} w_{i,v}.
\end{equation}
As evident from (\ref{eq:remains}), we will be interested in the
values of $\lambda^a_{j,v}$ only in the range $j =
n-k_a+1,\ldots,n$. We therefore rearrange these values by defining
a bijection $u: \cup_{a=1}^r \set{(a,n-k_a+1),\ldots,(a,n)}
\rightarrow \set{1,\ldots,l}$, and denote $\Lambda_{u(a,j),v} =
\lambda^a_{j,v}$. Plugging (\ref{eq:w}) into (\ref{eq:remains}),
we have:
\begin{eqnarray}
\nonumber \bigwedge_{a=1}^r \bigwedge_{i=d+1}^{n-l}
\bigwedge_{j=n-k_a+1}^{n} w^a_{i,j}  = \bigwedge_{i=d+1}^{n-l}
\bigwedge_{a=1}^r \bigwedge_{j=n-k_a+1}^{n} \sum_{v = n-l+1}^n \lambda^a_{j,v} w_{i,v} = \\
\nonumber \bigwedge_{i=d+1}^{n-l} \bigwedge_{u=1}^{l} \sum_{v =
n-l+1}^n \Lambda_{u,v} w_{i,v} = \bigwedge_{i=d+1}^{n-l}
det(\Lambda) w_{i,n-l+1} \wedge \ldots \wedge w_{i,n}.
\end{eqnarray}
The last transition is standard and is explained by the
skew-symmetry of the exterior product: all terms for which
$w_{i,v_1} \wedge \ldots \wedge w_{i,v_l}$ contains a recurring
$v_i = v_j$ are 0, and we are only left with the case $v_i =
\pi(i)$, where $\pi$ is a permutation of $\set{n-l+1,\ldots,n}$;
these terms are equal to $(-1)^{sign(\pi)} w_{i,n-l+1} \wedge
\ldots \wedge w_{i,n}$, producing the determinant of $\Lambda$.
Continuing, since $\Lambda$ does not depend on $i$ and using
(\ref{eq:basic-wedge}), we see that:
\[
\nonumber \bigwedge_{a=1}^r \bigwedge_{i=d+1}^{n-l}
\bigwedge_{j=n-k_a+1}^{n} w^a_{i,j}  = det(\Lambda)^{n-l-d}
\bigwedge_{i=d+1}^{n-l} \bigwedge_{j=n-l+1}^{n} w_{i,j} =
det(\Lambda)^{n-l-d} C^3_{n,l,d} d\mu^l_D(F),
\]
where $C^3_{n,l,d} = \abs{G(n-d,n-d-l)}$. To deduce
(\ref{eq:remains}), it remains to show that $\det(\Lambda) =
\Omega(\bar{E})$.

Recall that $\lambda^a_{j,v} = \scalar{e^a_j,g_v}$, and in the
range $j=n-k_a+1,\ldots,n$, these are exactly the coefficients of
the orthonormal bases $\bar{e^a} =
\set{e^a_{n-k_a+1},\ldots,e^a_n}$ of $E_a^\perp$ w.r.t. the
orthornormal basis $\bar{g} = \set{g_{n-l+1},\ldots,g_n}$ of
$F^\perp$. Using the orthogonality of $\bar{g}$, it is immediate
that $(\Lambda \Lambda^t)_{u(a_1,j_1),u(a_2,j_2)} =
\scalar{e^{a_1}_{j_1},e^{a_2}_{j_2}}$, and therefore
$det(\Lambda) = Vol_{F^\perp}(\bar{e})$ for $\bar{e} =
\set{\bar{e^1},\ldots,\bar{e^r}}$, which is exactly the
definition of $\Omega(\bar{E})$. Incidentally, this also shows
that $Vol_{F^\perp}(\bar{e})$ is invariant to taking an arbitrary
(not necessary orthonormal) basis $\bar{e^a}$ of $E_a^\perp$ with
$Vol_{E^\perp}(\bar{e^a})=1$, since this is easily checked for
$det(\Lambda)$.

\end{proof}

%%%%%%%%%%%%%%%%%%%%%%%%%%%%%%%%%%%%%%%%%%%%%%%%%%%%%%%%%%%%%%%%%%%%%%%%%%%%%%
%%%%%%%%%%%%%%%%%%%%%%%%%%%%%%%%%%%%%%%%%%%%%%%%%%%%%%%%%%%%%%%%%%%%%%%%%%%%%%

\bibliographystyle{amsalpha}
\bibliography{../../ConvexBib}

\end{document}